\documentclass{article}

\usepackage{amssymb,amscd,amsthm,verbatim,fancyhdr}
\usepackage{mathrsfs}
\usepackage{amsmath}

\begin{document}
\numberwithin{equation}{section}
\newtheorem{thm}{Theorem}
\newtheorem{lemma}{Lemma}
\newtheorem{clm}{Claim}
\newtheorem{remark}{Remark}
\newtheorem{definition}{Definition}
\newtheorem{cor}{Corollary}
\newtheorem{prop}{Proposition}
\newtheorem{statement}{Statement}


\newcommand{\hdt}{{\dot{\mathrm{H}}^{1/2}}}
\newcommand{\hdtr}{{\dot{\mathrm{H}}^{1/2}(\mathbb{R}^3)}}
\newcommand{\R}{\mathbb{R}}
\newcommand{\ei}{\mathrm{e}^{it\Delta}}
\newcommand{\ltrt}{{L^3(\mathbb{R}^3)}}
\newcommand{\ldrd}{{L^d(\mathbb{R}^d)}}
\newcommand{\lprd}{{L^p(\mathbb{R}^d)}}
\newcommand{\lt}{{L^3}}
\newcommand{\ld}{{L^d}}
\newcommand{\lp}{{L^p}}
\newcommand{\rt}{\mathbb{R}^3}
\newcommand{\rd}{\mathbb{R}^d}
\newcommand{\X}{\mathfrak{X}}
\newcommand{\F}{\mathfrak{F}}
\newcommand{\hdhalf}{{\dot H^\frac{1}{2}}}
\newcommand{\hdthalf}{{\dot H^\frac{3}{2}}}
\newcommand{\hdo}{\dot H^1}
\newcommand{\rthmiz}{\R^3\times(-\infty,0)}
\newcommand{\q}[2]{{#1}_{#2}}
\renewcommand{\t}{\theta}
\newcommand{\lxt}[2]{L_{x,\,t}^{#1}}
\newcommand{\rr}{\sqrt{x_1^2+x_2^2}}
\newcommand{\ve}{\varepsilon}
\newcommand{\hdhrt}{\dot H^\frac{1}{2}(\mathbb{R}^3)}
\renewcommand{\P}{\mathbb{P} }
\newcommand{\RR}{\mathcal{R} }
\newcommand{\TT}{\overline{T} }
\newcommand{\e}{\epsilon }
\newcommand{\D}{\Delta }
\renewcommand{\d}{\delta }
\renewcommand{\l}{\lambda }
\newcommand{\To}{\TT_1 }
\newcommand{\ukt}{u^{(KT)} }
\newcommand{\ttil}{\tilde{T_1} }
\newcommand{\etl}{e^{t\Delta} }
\newcommand{\et}{\mathscr{E}_T }
\newcommand{\se}{\mathscr{E}}
\newcommand{\ft}{\mathscr{F}_T }
\newcommand{\eti}{\mathscr{E}^{\infty}_T }
\newcommand{\fti}{\mathscr{F}^{\infty}_T }
\newcommand{\xjn}{x_{j,n}}
\newcommand{\xjpn}{x_{j',n}}
\newcommand{\ljn}{\l_{j,n} }
\newcommand{\ljpn}{\l_{j',n} }
\newcommand{\lkn}{\l_{k,n} }
\newcommand{\voj}{U_{0,j} }
\newcommand{\uoj}{U_{0,j} }
\newcommand{\voo}{U_{0,1} }
\newcommand{\uon}{u_{0,n} }
\newcommand{\N}{\mathbb{N} }
\newcommand{\Z}{\mathbb{Z} }
\newcommand{\E}{\mathscr{E} }
\newcommand{\tu}{\tilde{u} }
\newcommand{\tU}{\tilde{U} }
\newcommand{\etj}{{\E_{T^*_j}}}
\newcommand{\B}{\mathcal{B}}
\newcommand{\tujn}{\tilde{U}_{j,n}}
\newcommand{\tujpn}{\tilde{U}_{j',n}}
\newcommand{\soj}{\sum_{j=1}^J}
\newcommand{\soi}{\sum_{j=1}^\infty}
\newcommand{\hnj}{H_{n,J}}
\newcommand{\enj}{e_{n,J}}
\newcommand{\pnj}{p_{n,J}}
\newcommand{\lfoi}{L^5_{(0,\infty)}}
\newcommand{\lfhoi}{L^{{5/2}}_{(0,\infty)}}
\newcommand{\wnj}{w_n^J}
\newcommand{\rnj}{r_n^J}
\newcommand{\lfij}{{L^5_{I_j}}}
\newcommand{\lfiz}{{L^5_{I_0}}}
\newcommand{\lfio}{{L^5_{I_1}}}
\newcommand{\lfit}{{L^5_{I_2}}}
\newcommand{\lfhij}{{L^{{5/2}}_{I_j}}}
\newcommand{\lfhiz}{{L^{{5/2}}_{I_0}}}
\newcommand{\lfhio}{{L^{{5/2}}_{I_1}}}
\newcommand{\lfhit}{{L^{{5/2}}_{I_2}}}
\newcommand{\lfi}{{L^5_{I}}}
\newcommand{\lfhi}{{L^{{5/2}}_{I}}}
\newcommand{\doh}{D^\frac{1}{2}}
\newcommand{\tjn}{t_{j,n}}
\newcommand{\tjpn}{t_{j',n}}
\newcommand{\wnlj}{w_n^{l,J}}
\renewcommand{\O}{\mathcal{O}}
\newcommand{\bes}{{\dot B^{s_p}_{p,q}}}
\newcommand{\bespp}{{\dot B^{s_p}_{p,p}}}
\newcommand{\besa}{{\dot B^{s_{a}}_{a,b}}}
\newcommand{\besb}{{\dot B^{s_{p}}_{p,q}}}
\newcommand{\besr}{{\dot B^{s_{b}+\frac{2}{\rho}}_{b,r}}}
\newcommand{\tn}{{\tilde \|}}
\newcommand{\tp}{{\tilde{\phi}}}
\newcommand{\xto}{\xrightarrow[n\to\infty]{}}
\newcommand{\LL}{\Lambda}
\newcommand{\binf}{{\dot B^{-{d/p}}_{\infty,\infty}}}
\newcommand{\brq}{{\dot B^{s_{p,r}}_{r,q}}}
\newcommand{\bfin}{{\dot B^{s_{p,r}}_{r,q}}}
\newcommand{\trnj}{{ \tilde{r}_n^J}}
\def\longformule#1#2{
\displaylines{ \qquad{#1} \hfill\cr \hfill {#2} \qquad\cr } }
\def\sumetage#1#2{ \sum_{\binom{\scriptstyle {#1}}{\scriptstyle {#2}}}}
\def\longformule#1#2{ \displaylines{\qquad{#1} \hfill\cr \hfill {#2} \qquad\cr } }
\newcommand{\avint}[1]{{\int_{#1} \!\!\!\!\!\!\!\!\!\! -} \ \ }
\newcommand{\avintttt}[1]{{\int \!\!\!\!\!\!\! - \!\! \int_{#1} \!\!\!\!\!\!\!\!\!\!\!\!\!\!\!\!\!\!\!\!\!\!\!\! - } \ \ }
\newcommand{\avinttt}[1]{{\int \!\!\!\!\!\!\! - \!\! \int_{#1} \!\!\!\!\!\!\!\!\!\!\!\!\!\!\!\! - } \ \ }
\newcommand{\avintt}[1]{{\int \!\!\!\!\!\!\! - \!\! \int_{#1} \!\!\!\!\!\!\!\!\!\! - } \ \ }
\newcommand{\intt}[1]{{\int \!\!\! \int_{#1}}}
\newcommand{\g}{{\gamma}}
\newcommand{\n}{{\nabla}}
\newcommand{\p}{{\partial}}
\newcommand{\esssup}{\mathop{\mathrm{ess\,sup}}}

\title{Parabolic fractal dimension of forward-singularities for Navier-Stokes and liquid crystals inequalities}
\author{{\sc Gabriel S. Koch}\\ \small University of Sussex\\ \small Brighton BN1 9QH, UK \\ \small {\em g.koch@sussex.ac.uk}}

\date{}
\maketitle

\begin{abstract}
In 1985, Vladimir Scheffer discussed partial regularity results for what he called solutions to the ``Navier-Stokes inequality''.  These maps essentially satisfy the incompressibility condition as well as the local and global energy inequalities and the pressure equation which may be derived formally from the Navier-Stokes system of equations, but they are not required to satisfy the Navier-Stokes system itself.

In a previous work, the author extended this notion to a system considered by Fang-Hua Lin and Chun Liu in the mid 1990s related to models of the flow of nematic liquid crystals, which include the Navier-Stokes system when the `director field' $d$ is taken to be zero.  In addition to an extended Navier-Stokes system, the Lin-Liu model includes a further parabolic system which implies an a priori maximum principle for $d$ which they use to establish the same type of partial regularity result,  in terms of the parabolic Hausdorff dimension of sets, as that which was known in the Navier-Stokes setting.  For the analogous `inequality' one loses this maximum principle and the author previously explored the consequences for such Hausdorff-dimensional partial regularity.

The current work explores similar consequences for partial regularity, but with respect to the `parabolic fractal dimension' $\textrm{dim}_{\textrm{pf}}$ (also called the upper box-counting, capacity or Minkowski dimension). In 2018, relying (as did Lin and Liu) on the boundedness of $d$ coming from the maximum principle, Qiao Liu proved that solutions to the Lin-Liu model satisfy  ${\textrm{dim}_{\textrm{pf}}(\Sigma_{-} \cap \mathcal{K}) \leq  \tfrac{95}{63}}$ for any compact set $\mathcal{K}$, where  $\Sigma_{-}$ is the set of `forward-singular' space-time points near which the solution blows up forwards in time.  For solutions to the corresponding `inequality', the author proves here that, without any compensation for the lack of maximum principle, one has $\textrm{dim}_{\textrm{pf}}(\Sigma_{-} \cap \mathcal{K}) \leq  \tfrac {55}{13}$.  A range of  criteria is also established, including as just one example the boundedness of $d$, any one of which is shown to furthermore imply  that  $\textrm{dim}_{\textrm{pf}}(\Sigma_{-} \cap \mathcal{K}) \leq  \tfrac{95}{63}$ for solutions to the inequality, just as Q. Liu proved for solutions to the full Lin-Liu system.
\end{abstract}

\section{Introduction}
\noindent
In \cite{linliu95} and \cite{linliu}, Fang-Hua Lin and Chun Liu prove existence and partial regularity (similar to results in \cite{caf} for Navier-Stokes) of certain solutions to the following system, which reduces to the classical Navier-Stokes system in the case $d\equiv 0$ (here we have set various parameters equal to one for simplicity):
\begin{equation}\label{maineq}
\!\!\boxed{\begin{array}{rcl}
u_t  - \Delta u + \nabla^T \cdot [u\otimes u+\n d \odot \n d ]+ \nabla p & = & 0\\\\
\nabla \cdot u & = & 0\\\\
d_t -\Delta d + (u\cdot \nabla)d  +f(d) &=&0
\end{array}}
\end{equation}
with $f=\n F$ for the scalar field $F$ given by $F(x) := (|x|^2-1)^2$
so that
$f(x)= 4(|x|^2-1)x$
(and in particular $f(0) = 0$).  In \eqref{maineq}, for vector fields $v$ and $w$, the matrix fields $v\otimes w$ and $\n v \odot \n w$ are defined to be the ones with entries
$$(v\otimes w)_{ij} = v_i w_j \quad \textrm{and} \quad (\n v \odot \n w)_{ij}=v_{,i}\cdot w_{,j}:= \frac{\partial v_k}{\partial x_i}\frac{\partial w_k}{\partial x_j}$$
(summing over the repeated index $k$ as per the Einstein convention), and for a matrix field $J= (J_{ij})$, we define
the vector field $\n^T \cdot J$ by
$$(\n^T \cdot J)_i:=J_{ij,j}:= \frac {\partial J_{ij}}{\partial x_j}$$
(summing again over $j$), where we think formally of $\n$ (as well as any vector field) as a column vector and $\n^T$ as a row vector.
\\\\
We take the spatial dimension to be three, so that for some $\Omega \subseteq \R^3$ and $T>0$, we are seeking maps of the form
$$u,d : \Omega \times (0,T) \to \R^3 \quad \textrm{and} \quad \quad p: \Omega \times (0,T) \to \R$$
satisfying \eqref{maineq}, where
$$F:\R^3 \to \R  \quad \textrm{and} \quad f:\R^3 \to \R^3$$
are fixed as above.  As usual, $u$ represents the velocity vector field of a fluid, $p$ is the scalar pressure in the fluid, and, as in nematic liquid crystals models, $d$ corresponds roughly to the `director field' representing the local orientation of rod-like molecules, with $u$ also giving the velocities of the centers of mass of those anisotropic molecules.
\ \\\\
We note that by formally taking the divergence $\n \cdot$ of the first line in (\ref{maineq}) we obtain the usual `pressure equation'
\begin{equation}\label{preseqa}
-\D p = \n \cdot (\n^T \cdot [u\otimes u+\n d \odot \n d ])\, .
\end{equation}
As in the Navier-Stokes ($d\equiv 0$) setting, one may formally deduce (see e.g. \cite{koch2021} for more details) from (\ref{maineq}) the following global and local energy inequalities which one may expect `sufficiently nice' solutions of (\ref{maineq}) to satisfy:
\begin{equation}\label{globenineq}
\displaystyle{\frac{d}{d t}\int_{\Omega} \left[ \frac{|u|^2}2 +\frac{|\n d|^2}2 + F(d)  \right]\, dx + \int_{\Omega}\left[|\n u|^2 + |\D d-f(d)|^2 \right]}\, dx\leq 0
\end{equation}
for each $t\in (0,T)$, as well as a localized version
\begin{equation}\label{locenineq}
\begin{array}{l}
\displaystyle{\frac{d}{d t}\int_{\Omega}\left[\left(\frac{|u|^2}{2} + \frac{|\n d|^2}2 \right)\phi \right] \, dx + \int_{\Omega}\left(|\n u|^2+ |\n^2 d|^2\right)\phi \, dx}
\\\\
\qquad \qquad \qquad
\begin{array}{l}
\displaystyle{\leq  \int_{\Omega} \bigg[\left(\frac{|u|^2}{2}+ \frac{|\n d|^2}{2}\right)(\phi_t + \D \phi) +\left(\frac{|u|^2}2 +\frac{|\n d|^2}2 + p\right)u\cdot \n \phi} \\\\
\displaystyle{ \qquad  \qquad  \qquad \qquad \qquad +\ \  u\otimes \n \phi : \n d \odot \n d  \ \  -
\ \ \phi\n^T [f(d)]:  \n^T d\bigg]\, dx}
\end{array}
\end{array}
\end{equation}
for $t\in (0,T)$ and each smooth, compactly supported in $\Omega$ and non-negative scalar field $\phi \geq 0$.  (For Navier-Stokes, i.e. when $d\equiv 0$, one may omit all terms involving $d$, even though $0\neq F(0)\notin L^1(\R^3)$.)  As explained in \cite{koch2021}, a Gr\"onwall argument shows that this also formally implies a local energy inequality of the form (\ref{locenta}) below for some $\bar C>0$ depending only on $T$.
\\\\
In the Navier-Stokes setting, it was asserted by Vladimir Scheffer in \cite{scheffer3} that in fact the proof of the partial regularity result (see below) of Caffarelli-Kohn-Nirenberg \cite{caf}  does not require the full set of equations in (\ref{maineq}).  He mentions that the key ingredients are membership of the global energy spaces, the local energy inequality (\ref{locenineq}), the divergence-free condition $\n \cdot u =0$ and the {\em pressure} equation (\ref{preseqa}) (with $d\equiv 0$ throughout).  Scheffer called pairs $(u,p)$ satisfying these four requirements solutions to the ``Navier-Stokes inequality'', equivalent to solutions to the Navier-Stokes equations with a forcing $g$ which satisfies $g\cdot u \leq 0$ everywhere.  (The Navier-Stokes solutions considered in \cite{caf}, called ``suitable weak solutions'', are roughly speaking solutions to Scheffer's  `inequality' which moreover satisfy \eqref{maineq} itself -- with, of course, $d\equiv 0$.)  In contrast, the results in \cite{linliu} (for a similar class of `suitable weak solutions' to \eqref{maineq} for more general $d$) do very strongly use the third equation in (\ref{maineq}) in that it implies a maximum principle for $d$, making it natural to assume $d\in L^\infty(\Omega \times (0,T))$.
\\\\
In this paper, as in the author's previous paper \cite{koch2021}, we continue to explore what happens if one considers the analog of Scheffer's ``Navier-Stokes inequality''  for the system (\ref{maineq}) when $d\neq 0$.  That is, we consider triples $(u,d,p)$ with global regularities implied (at least when $\Omega$ is bounded and under suitable assumptions on the initial data) by (\ref{globenineq}) which satisfy (\ref{preseqa}) and $\n \cdot u=0$ weakly as well as (a formal consequence of) (\ref{locenineq}), but are {\em not} necessarily weak solutions of the first and third equations (i.e., the two vector equations) in (\ref{maineq}).  In particular, we will {\em not} assume that $d\in L^\infty(\Omega \times (0,T))$, which would have been reasonable in view of the third equation in (\ref{maineq}).  Specifically,  we will address the following scenario:
\\\\
Fixing an open set $\Omega \subseteq \R^3$ and $T,\bar C\in (0,\infty)$ and setting $\Omega_T:=\Omega \times (0,T)$, we consider
$u,d:\Omega_T  \to \R^3$ and $p: \Omega_T \to \R$ which satisfy the following four assumptions\footnote{See \cite{koch2021} for more explanation as to why these assumptions are natural from the viewpoint of \eqref{maineq}.}.
\begin{enumerate}
\item $u$, $d$ and $p$  belong to the following spaces:\footnote{For a vector field $f$ or matrix field $J$ and scalar function space $X$, by $f\in X$ or $J\in X$ we mean that all components or entries of $f$ or $J$ belong to $X$; by $\n^2 f \in X$ we mean all second partial derivatives of all components of $f$ belong to $X$; etc.}
\begin{equation}\label{enspaces}
d\in L^\infty(0,T;L^4(\Omega))\,  , \quad
u, \n d \in L^\infty(0,T;L^2(\Omega))\, , \quad
\nabla u,  \n^2 d  \in L^2(\Omega_T)
\end{equation}
and
\begin{equation}\label{pspace}
p\in L^{\frac 53}(\Omega_T)\, ;
\end{equation}
\item $u$ is weakly\footnote{Locally integrable functions, e.g. $u \in  L^1_{\mathrm{loc}}(\Omega_T)$ in \eqref{divfree} and $u\otimes u + \n d \odot \n d \in  L^1_{\mathrm{loc}}(\Omega_T)$ in \eqref{preseq},  will always be associated to the standard distribution whose action is integration against a suitable test function so that, e.g., $[\n \cdot u](\psi)= -[u](\n \psi):= -\int u\cdot \n \psi$ for $\psi \in \mathcal{D}(\Omega_T)$. }  divergence-free:
\begin{equation}\label{divfree}
\nabla \cdot u = 0 \quad \textrm{in} \quad \mathcal{D}'(\Omega_T)\, ;
\end{equation}
\item The following pressure equation holds weakly:
\begin{equation}\label{preseq}
-\Delta p = \nabla \cdot [\nabla^T \cdot (u\otimes u+\n d \odot \n d)]\quad \textrm{in} \quad \mathcal{D}'(\Omega_T)\, ;
\end{equation}
\item The following  local energy inequality holds:\footnote{For brevity, for $\omega \subseteq \R^3$, we set  $$\int_{\omega \times \{t\}} g\, dx:=\int_{\omega}g(x, t)\, dx\, .$$ \label{fnintdef}}
\begin{equation}\label{locenta}
 \boxed{\begin{array}{l}
 \int_{\Omega \times \{t\}} \left(|u|^2 + |\n d|^2\right) \phi \, dx + \int_{0}^t \int_{\Omega} \left(|\n u|^2+ |\n^2 d|^2\right) \phi \, dx\, d\tau  \\\\
\qquad \leq \bar C\int_{0}^t \big\{ \int_{\Omega \times \{\tau\}} \left[\left(|u|^2 + |\n d|^2\right)|\phi_t + \D \phi|  + |d|^2 |\n d|^2\phi \right]\, dx\\\\
 \qquad \qquad \quad +\  \  \big|\int_{\Omega \times \{\tau\}} \big[\big(\frac{|u|^2}2 +\frac{|\n d|^2}2 + p\big)u\cdot \n \phi + u\otimes \n \phi : \n d \odot \n d\big]\, dx\big|\ \big\}\, d\tau\ \\\\
\qquad \textrm{for}\ \textrm{a.e.}\ t\in (0,T)\  \quad \textrm{and}\quad \forall \ \phi \in  \mathcal{C}_0^\infty(\Omega \times (0,\infty))\ \textrm{s.t.}\ \phi \geq 0\, .
\end{array}}
\end{equation}
\end{enumerate}

\begin{remark}\label{pspacesremark}
In the case $\Omega = \rt$, the condition (\ref{pspace}) on the pressure follows from (\ref{enspaces}) and (\ref{preseq}) if $p$ is taken to be the potential-theoretic solution to (\ref{preseq}), since (\ref{enspaces}) implies that  ${u,\n d\in L^{\frac{10}3}(\Omega_T)}$ by interpolation  and Sobolev embeddings, and then (\ref{preseq}) gives   $p\in L^{\frac{5}3}(\Omega_T)$ by Calderon-Zygmund estimates.  For a more general $\Omega$, the existence of such a $p$ can be derived from the motivating equation (\ref{maineq})  (e.g. by estimates for the Stokes operator), see \cite{linliu} and the references therein.  Here, however, we will not refer to (\ref{maineq}) at all and simply {\em assume} $p$ satisfies (\ref{pspace}) and address the partial regularity of such a hypothetical set of functions satisfying  (\ref{enspaces}) - (\ref{locenta}).
\end{remark}
\noindent
The partial regularity results of \cite{caf,linliu,koch2021} are given in terms of the (outer) parabolic Hausdorff measure $\mathcal{P}^k$ defined  as follows (see \cite[pp.783-784]{caf}), where $B_r(x)$ is the standard Euclidean ball in $\R^3$ of radius $r>0$ centered at $x\in \R^3$:
\begin{definition}[Parabolic Hausdorff measure and dimension]\label{phausdorff}
For any set $\mathcal{A} \subset \R^3 \times \R$ and $k\geq 0$, its $k$-dimensional (outer) parabolic Hausdorff measure $\mathcal{P}^k(\mathcal{A})$ is defined as
$$\mathcal{P}^k(\mathcal{A}):=\lim_{\delta \searrow 0}\mathcal{P}^k_\delta(\mathcal{A})\, ,$$
where
$$\mathcal{P}^k_\delta(\mathcal{A}):=\inf\left\{\, \sum_{j=1}^\infty r_j^k\ \bigg| \ \mathcal{A}\subset \bigcup_{j=1}^\infty Q_{r_j}\, , r_j < \delta\  \forall j\in \N\, \right\}$$
and $Q_r$ is any parabolic cylinder of `radius' $r>0$, i.e.
\begin{equation}\label{standcyldef}
Q_r=Q_r(x,t):= B_r(x) \times (t-r^2,t) \subset \R^3 \times \R
\end{equation}
 for some $x\in \R^3$ and $t\in \R$.
The parabolic Hausdorff dimension $\textrm{dim}_{\mathcal{P}}(\mathcal{A})$ of $\mathcal{A}$ is then defined as
$$\textrm{dim}_{\mathcal{P}}(\mathcal{A}):=\inf \{\, k\geq 0 \, | \, \mathcal{P}^k(\mathcal{A})=0\, \}\, .$$
\end{definition}
\noindent
The results in \cite{caf} and \cite{linliu} state roughly that for suitable weak solutions to \eqref{maineq} (including the Navier-Stokes setting $d\equiv 0$ of \cite{caf}), the singular set $\Sigma$ of space-time points $(x,t)$ about each of which $u$ and $\n d$ are unbounded in some space-time neighborhood satisfies $\mathcal{P}^1(\Sigma)=0$ (hence $\textrm{dim}_{\mathcal{P}}(\Sigma)\leq 1$).  In these settings the vector-field $d$ is bounded due to \eqref{maineq}, while in \cite{koch2021}, for solutions
to the corresponding `inequality', the author proved that (without any compensation for the lack of boundedness of $d$) one has $\mathcal{P}^{\frac 92 + \d}(\Sigma)=0$ for any $\d >0$ (hence $\textrm{dim}_{\mathcal{P}}(\Sigma)\leq \tfrac 92$); the author moreover provided a range of  criteria, including as just one example the boundedness of $d$, any one of which  was shown to furthermore imply (as in \cite{caf,linliu}) that  $\mathcal{P}^1(\Sigma)=0$ even for solutions to the inequality.
\ \\\\
In \cite{qliu2018}, on the other hand, Qiao Liu  proves (along the lines of other papers such as, e.g., \cite{kukpei2012} and the references therein), for the same type of solutions to \eqref{maineq} considered by \cite{linliu}, a partial regularity result in terms of the parabolic fractal dimension (also called the upper box-counting or capacity dimension, see \cite{kukpei2012} and the references therein, as well as the upper Minkowski dimension as it is referred to in \cite{qliu2018}) defined as follows:
\begin{definition}[Parabolic fractal dimension\footnote{It is pointed out in \cite{qliu2018} that, in general, $\textrm{dim}_{\mathcal{P}}(\mathcal{A}) \leq \textrm{dim}_{\textrm{pf}}(\mathcal{A})$, where strict inequality is known to occur in some instances. }, see e.g. \cite{kukpei2012}]\label{parfracdimdef}
For any bounded set $\mathcal{A}\subseteq \R^3 \times \R$, its parabolic fractal dimension $\textrm{dim}_{\textrm{pf}}(\mathcal{A})$  is defined as
$$\textrm{dim}_{\textrm{pf}}(\mathcal{A}):=\limsup_{r\searrow 0} \frac{\log N(\mathcal{A},r)}{\log (\tfrac 1r)}$$
where, for any $r>0$,
$$N(\mathcal{A},r):=\min \left\{N\in \N \ \bigg| \ \mathcal{A} \subseteq \bigcup_{j=1}^N Q_r^*(z_j) \quad \textrm{for some} \ \{z_j\}_{j=1}^N \subset \R^3 \times \R \right\}$$
is the minimum number of (centered) parabolic cylinders of the form\footnote{Note that in view of \eqref{standcyldef} and \eqref{qstarsdef} one always has $Q_r(z) \subset Q_r^*(z)$ for any $z=(x,t)\in \R^3 \times \R$ and $r>0$.}
\begin{equation}\label{qstarsdef}
Q^{*}_{r}(x,t):=B_r(x)\times (t- r^2,t+r^2) \qquad (x\in \R^3, t\in \R, r>0)
\end{equation}
with the fixed `radius' $r$ required to cover $\mathcal{A}$.
\end{definition}
\noindent
Roughly speaking, for suitable weak solutions to \eqref{maineq}, Q. Liu proves (see Theorem \ref{corbliuthm} below) that $\textrm{dim}_{\textrm{pf}}(\Sigma_{-} \cap \mathcal{K}) \leq  \tfrac{95}{63}$ for any compact space-time set $\mathcal{K}$, where  $\Sigma_{-}$ is the set of `forward-singular' space-time points, near each of which the solution blows up forwards in time.  As we will describe below, the proof of Theorem \ref{corbliuthm} given in \cite{qliu2018} relies again on the fact that $d$ is bounded, a consequence of the third equation in \eqref{maineq}. In contrast, for solutions to the corresponding `inequality', we will prove in Theorem \ref{corbkochthm} below  that, without any compensation for the lack of maximum principle and hence when $d$ is potentially unbounded, one can at least prove $\textrm{dim}_{\textrm{pf}}(\Sigma_{-} \cap \mathcal{K}) \leq  \tfrac {55}{13}$. It is reasonable that Theorem \ref{corbkochthm}, which may be viewed as an adaptation of Q. Liu's  Theorem \ref{corbliuthm}, yields a weaker conclusion (the bound $\tfrac {55}{13}$ rather than the smaller $\tfrac{95}{63}$), as one has removed the  boundedness of $d$ from the set of essential properties used in the proof.  As in \cite{koch2021}, we will also establish (in Theorem \ref{corbkochbthmb} below) a range of  criteria (namely \eqref{gsigmadeffina}) for solutions to the inequality, including as just one example the boundedness of $d$, any one of which  would furthermore imply (without any reference to \eqref{maineq} itself) that  $\textrm{dim}_{\textrm{pf}}(\Sigma_{-} \cap \mathcal{K}) \leq  \tfrac{95}{63}$; Theorem \ref{corbkochbthmb} may therefore be viewed as a generalization of Q. Liu's \linebreak Theorem \ref{corbliuthm}.
\\\\
Our key observation  which allows us to work without any maximum principle is that,  in view of the global energy (\ref{globenineq}) and the particular forms of $F$ and $f$, it is reasonable (see \cite{koch2021}) to assume (\ref{enspaces}) which implies (at least for bounded $\Omega$ and  $T<\infty$) that ${d \in L^\infty(0,T;L^6(\Omega))} \cap L^{10}(\Omega \times (0,T))$,  which is sufficient for our purposes.   Indeed, the assumptions in (\ref{enspaces}) imply (locally in space, or if $\Omega$ is bounded) that $d \in L^\infty(0,T;H^1(\Omega)) \hookrightarrow L^\infty(0,T;L^6(\Omega))$ by Sobolev embedding, as well as the standard fact (due to interpolation and Sobolev embedding, see e.g. \cite{koch2021} for details) that $\n d \in L^{\frac 2\alpha}(0,T;L^{\frac 6{3-2\alpha}}(\Omega))$ for any $\alpha \in [0,1]$ since $\n d$ lives in the same Navier-Stokes-type energy spaces as does $u$.  (Taking $\alpha:=\frac 35$ gives the usual fact that $u,\n d \in L^{\frac{10}3}(\Omega \times (0,T))$.) Taking $\alpha:=\frac 15$, we see that (for bounded $\Omega$ and $T<\infty$) $d\in L^{10}(0,T;W^{1,\frac{30}{13}}(\Omega)) \hookrightarrow L^{10}(\Omega \times (0,T))$ by Sobolev embedding.

\section{Main results and context}
\noindent
In \cite{qliu2018}, Q. Liu proves (for the spatial setting $\Omega = \R^3$) a theorem (\cite[Theorem 1]{qliu2018}) which may  be stated essentially as follows:
\begin{thm}[\cite{qliu2018}, Theorem 1]\label{corbliuthm}
Fix an open set\footnote{As mentioned above, in  \cite{qliu2018} it is actually assumed that $\Omega = \R^3$, but it is clear from the proof that one may  establish  Theorem \ref{corbliuthm} as stated here; this is in fact carried out in a more general setting in Theorem \ref{corbkochbthmb} below.} $\Omega \subseteq \R^3$ and $T, \bar C \in (0,\infty)$, set  $\Omega_T:=\Omega \times (0,T)$, suppose
$u,d:\Omega_T  \to \R^3$ and $p: \Omega_T \to \R$ satisfy the assumptions \eqref{enspaces}, \eqref{pspace}, \eqref{divfree},  \eqref{preseq} and \eqref{locenta} and suppose additionally that $d\in L^\infty(\Omega_T)$
(and that \eqref{maineq} holds\footnote{In fact, \eqref{maineq}  is not needed, see Theorem \ref{corbkochbthmb} below.} weakly).
Define the forward(-in-time)-singular set\footnote{Recall the definition \eqref{standcyldef} used to define the (standard parabolic) cylinder $Q_r(z_0)$ which appears in the definition of $\Sigma_{-}$.  In the notation $\Sigma_{-}$, the ``$-$'' can therefore be viewed as indicating blowup for times approaching $t_0 \in (0,T)$ from the left, i.e. moving forward in time, for $z_0 = (x_0,t_0) \in \Sigma_{-}$.  } $\Sigma_{-} \subseteq \Omega_T$ by
$$\Sigma_{-}:= \{\ z_0 \in \Omega_T \ | \ \left\|\,  |u| + |\n d| \, \right\|_{L^\infty(Q_r(z_0)\cap \Omega_T)} = +\infty \ \ \forall r>0\ \}\, .$$
Then for any compact set $\mathcal{K} \subset \Omega_T$,
$$\textrm{dim}_{\textrm{pf}}(\Sigma_{-} \cap \mathcal{K}) \leq \tfrac 53 - \gamma \quad \textrm{for any}\ \gamma \in (0,\tfrac {10}{63})\, ,$$
so that\footnote{During the preparation of this paper, the author discovered that the bound in \eqref{partialregtyp} appears to have been recently improved slightly in \cite{qliu2021} from $\textrm{dim}_{\textrm{pf}}(\Sigma_{-} \cap \mathcal{K}) \leq  \tfrac{95}{63} \approx 1.51$ to  $\textrm{dim}_{\textrm{pf}}(\Sigma_{-} \cap \mathcal{K}) \leq  \frac{835}{613} \approx 1.36$.}
\begin{equation}\label{partialregtyp}
\textrm{dim}_{\textrm{pf}}(\Sigma_{-} \cap \mathcal{K}) \leq \tfrac 53 - \tfrac {10}{63}= \tfrac{95}{63}\, .
\end{equation}
\end{thm}
\noindent
Here, we will prove two related theorems.  The first one (Theorem \ref{corbkochbthmb}) shows that one does not need to assume that \eqref{maineq} itself holds as long as one assumes that $d$ is either bounded (as in Theorem \ref{corbliuthm}) or satisfies an alternative assumption (see \eqref{gsigmadeffina}) which plays the same role in the proof:
\begin{thm}\label{corbkochbthmb}
Fix an open set $\Omega \subseteq \R^3$ and $T, \bar C \in (0,\infty)$,  set $\Omega_T:=\Omega \times (0,T)$ and suppose \linebreak
$u,d:\Omega_T  \to \R^3$ and $p: \Omega_T \to \R$ satisfy the assumptions \eqref{enspaces}, \eqref{pspace}, \eqref{divfree},  \eqref{preseq} and \eqref{locenta}.
Suppose as well that
\begin{equation}\label{gsigmadeffina}
g_\sigma:=\sup_{\{r,z_0\, | \, Q_r(z_0)\subseteq \Omega_T\}}\left(  \frac 1{r^{2+\frac \sigma2}}
 \intt{Q_r(z_0)}|d|^\sigma |\n d|^{3(1-\frac \sigma 6)}\, dz\right)<\infty \quad \textrm{for some} \ \sigma \in (5,6]\, .
\end{equation}
Define the forward-singular set $\Sigma_{-} \subseteq \Omega_T$ by
$$\Sigma_{-}:= \{\ z_0 \in \Omega_T \ | \ \left\|\,  |u| + |\n d| \, \right\|_{L^\infty(Q_r(z_0)\cap \Omega_T)} = +\infty \ \ \forall r>0\ \}\, .$$
Then for any compact set $\mathcal{K} \subset \Omega_T$,
$$\textrm{dim}_{\textrm{pf}}(\Sigma_{-} \cap \mathcal{K}) \leq \tfrac 53 - \gamma \quad \textrm{for any}\ \gamma \in (0,\tfrac {10}{63})\, ,$$
so that
$$\textrm{dim}_{\textrm{pf}}(\Sigma_{-} \cap \mathcal{K}) \leq \tfrac{95}{63}\, .$$
\end{thm}
\noindent
Note that if $d\in L^\infty(\Omega_T)$ (as in Theorem \ref{corbliuthm}), then\footnote{In fact, the two assumptions are equivalent for $\sigma=6$ due to Lebesgue's differentiation theorem.} \eqref{gsigmadeffina} holds with $\sigma=6$  and hence Theorem \ref{corbkochbthmb} in particular implies Theorem \ref{corbliuthm}.  The proofs of Theorem \ref{corbliuthm} and Theorem \ref{corbkochbthmb} are, however, very similar, and the endpoint $\sigma=6$ is included in \eqref{gsigmadeffina} mainly to point out that the proof of Theorem \ref{corbliuthm} given by Q. Liu in \cite{qliu2018} does not in fact require the full set of assumptions (as will be essentially clear from the proof of Theorem \ref{corbkochbthmb} given below).
\\\\
The second theorem (Theorem \ref{corbkochthm}) addresses the scenario where one removes the assumption \eqref{gsigmadeffina} altogether from Theorem \ref{corbkochbthmb} (or, indeed, from Theorem \ref{corbliuthm} when $\sigma=6$):
\begin{thm}\label{corbkochthm}
Fix an open set $\Omega \subseteq \R^3$ and $T, \bar C \in (0,\infty)$,  set $\Omega_T:=\Omega \times (0,T)$ and suppose \linebreak
$u,d:\Omega_T  \to \R^3$ and $p: \Omega_T \to \R$ satisfy the assumptions \eqref{enspaces}, \eqref{pspace}, \eqref{divfree},  \eqref{preseq} and \eqref{locenta}.
Define the forward-singular set
 $\Sigma_{-} \subseteq \Omega_T$ by
$$\Sigma_{-}:= \{\ z_0 \in \Omega_T \ | \ \left\|\,  |u| + |\n d| \, \right\|_{L^\infty(Q_r(z_0)\cap \Omega_T)} = +\infty \ \ \forall r>0\ \}\, .$$
Then for any compact set $\mathcal{K} \subset \Omega_T$,
$$\textrm{dim}_{\textrm{pf}}(\Sigma_{-} \cap \mathcal{K}) \leq 5-\delta \quad \textrm{for any}\ \delta \in (\tfrac 12,\tfrac {10}{13})\, ,$$
so that
$$\textrm{dim}_{\textrm{pf}}(\Sigma_{-} \cap \mathcal{K}) \leq  5-\tfrac {10}{13} = \tfrac {55}{13}\, .$$
\end{thm}
\noindent
The main mechanism used to establish results such as Theorems \ref{corbliuthm}, \ref{corbkochbthmb} and \ref{corbkochthm} is
the following general proposition which was already available in the literature  prior to \cite{qliu2018} (see, e.g., \cite{kukpei2012} and the references therein):
\begin{prop}\label{covarglem}
Fix any open set $\Omega \subseteq \R^3$ and $\lambda, \bar r,c_0,T\in (0,\infty)$ and set $\Omega_T:=\Omega \times (0,T)$. Suppose $\mathcal{S} \subseteq \Omega_T$  and that $H:\Omega_T \to [0,\infty]$ is a non-negative Lebesgue-measurable function such that
\begin{equation}\label{partialreggenfb}
0\leq H\in L^1(\Omega_T)\, ,
\end{equation}
and suppose that the following property holds in general (recall (\ref{qstarsdef})):
\begin{equation}\label{smalllam}
\left .\begin{array}{c}
z_0\in \mathcal{S}\\
0< r \leq \bar r\\
Q^*_{r}(z_0) \subseteq \Omega_T
\end{array} \right\} \quad \Longrightarrow \quad \frac 1{r^{\lambda}}\int_{Q^*_{r}(z_0)} H(z)\, dz \geq   c_0\, .
\end{equation}
Then for any compact set $\mathcal{K} \subset \Omega_T$,
$$\textrm{dim}_{\textrm{pf}}(\mathcal{S}\cap \mathcal{K}) \leq \lambda\, .$$
\end{prop}
\noindent
(A proof of Proposition \ref{covarglem} is roughly outlined in \cite{qliu2018}; for completeness we will prove Proposition \ref{covarglem} carefully in Section \ref{proponepfsec} below.)
\\\\
The novel and essential elements of Theorems \ref{corbliuthm}, \ref{corbkochbthmb} and \ref{corbkochthm} are therefore the following more specific `epsilon-regularity' type lemmas, each of which fairly immediately implies the corresponding theorem in view of Proposition  \ref{covarglem}.  The first, proved in \cite{qliu2018}, is the essence of \cite[Theorem 1]{qliu2018} (i.e., Theorem \ref{corbliuthm} above) and may be stated essentially as follows:

\begin{lemma}[\cite{qliu2018}, Lemma 6]\label{corbliu}
Fix any $\gamma \in (0,\tfrac {10}{63})$ and $\bar C, D \in (0,\infty)$.  There exist numbers  ${\e^*=\e^*(\gamma,\bar C,D) \in (0,1)}$ and  $r^*=r^*(\gamma) \in (0,1)$ so small
that the following holds for any fixed open set\footnote{As noted above, Q. Liu \cite{qliu2018} actually assumes $\Omega = \R^3$ so that the Calderon-Zygmund estimates apply directly to \eqref{preseq}; however this is unnecessary (see Lemma \ref{corbkochb} below) in the present context in view of the assumption \eqref{pspace} on the pressure.} $\Omega \subseteq \R^3$ and $T\in (0,\infty)$:
\\\\
Set  $\Omega_T:=\Omega \times (0,T)$, suppose
$u,d:\Omega_T  \to \R^3$ and $p: \Omega_T \to \R$ satisfy the assumptions \eqref{enspaces}, \eqref{pspace}, \eqref{divfree},  \eqref{preseq} and \eqref{locenta} and suppose additionally that $$d\in L^\infty(\Omega_T) \quad \textrm{with} \quad \|d\|_{L^\infty(\Omega_T)} \leq D$$
(and that \eqref{maineq} holds\footnote{In fact, \eqref{maineq}  is not needed, as in  Lemma \ref{corbkochb} below.} weakly).  For any $z_0\in \Omega_T$ and ${r_1} \in (0,r^*]$ such that $Q_{r_1}(z_0) \subseteq \Omega_T$, if
$$
 \frac 1{r_1^{\frac 53 - \gamma }} \int_{Q_{r_1}(z_0)} \left(|\n u|^2 + |\n^2 d|^2 + |u|^{\frac {10}3} +|\n d|^{\frac {10}3} + |p|^{\frac {5}3}\right) \, dz \leq   \e^*\, ,
$$
then
$$u, \n d \in L^\infty(Q_{\frac {r_0}2}(z_0))
$$
for some small $r_0\in (0,r_1)$.
\end{lemma}
\noindent
\ \\
Similarly, the essence of Theorems  \ref{corbkochbthmb} and \ref{corbkochthm} are the following Lemmas \ref{corbkochb} and \ref{corbkoch} which we will prove in Section  \ref{lemmasproofssec} below:
\begin{lemma}\label{corbkochb}
Fix any $\gamma \in (0,\tfrac {10}{63})$, $\sigma \in (5,6]$ and $\bar C, D \in (0,\infty)$. There exist numbers \linebreak $\e^*=\e^*(\gamma,\sigma,\bar C, D) \in (0,1)$ and  $r^*=r^*(\gamma) \in (0,1)$ so small
that the following holds for any fixed open set $\Omega \subseteq \R^3$ and $T\in (0,\infty)$:
\\\\
Set $\Omega_T:=\Omega \times (0,T)$, suppose
$u,d:\Omega_T  \to \R^3$ and $p: \Omega_T \to \R$ satisfy the assumptions \eqref{enspaces}, \eqref{pspace}, \eqref{divfree},  \eqref{preseq} and \eqref{locenta} and suppose additionally  that
\begin{equation}\label{gsigmadeffin}
g_\sigma  <\infty \quad \textrm{with} \quad g_\sigma \leq D\, ,
\end{equation}
with $g_\sigma$ defined as in (\ref{gsigmadeffina}).  For any $z_0\in \Omega_T$ and ${r_1} \in (0,r^*]$ such that $Q_{r_1}(z_0) \subseteq \Omega_T$, 
if\, \!\!\footnote{Note, in contrast to Lemma \ref{corbliu}, the appearance of the additional term $|d|^{10}$ in the smallness assumption in both Lemma \ref{corbkochb} and Lemma \ref{corbkoch}.  This is the key new idea which allows us to remove the requirement that $d\in L^\infty$.}
$$
 \frac 1{r_1^{\frac 53 - \gamma }} \int_{Q_{r_1}(z_0)} \left(|\n u|^2 + |\n^2 d|^2 + |u|^{\frac {10}3} +|\n d|^{\frac {10}3} + |p|^{\frac {5}3} + |d|^{10}\right) \, dz \leq   \e^*\, ,
$$
then
$$u, \n d \in L^\infty(Q_{\frac {r_0}2}(z_0))
$$
for some small $r_0\in (0,r_1)$.
\end{lemma}
\noindent
\ \\
Note that, in the particular case when $d\in L^\infty(\Omega_T)$, one can take $\sigma:=6$ and $D:=|Q_1|\|d\|_{L^\infty(\Omega_T)}$ in Lemma \ref{corbkochb} to recover the type of result in Lemma \ref{corbliu}.  \\\\\\



\begin{remark}\label{remthree}
For $\g \approx \frac{10}{63}$, setting
$$ M_{\g}:= \frac{154}{27(1-\frac {63}{10}\g)}
\quad \textrm{and} \quad N\in \N \cap [M_{\g},M_{ \g}+1)\, ,$$
it will be clear from the proof of Lemma \ref{corbkochb} that its conclusion is true, for example, with
$$
r^* := (\tfrac 14)^{6M_{ \g}+11} \quad \textrm{and} \quad r_0:=
2(\tfrac {r_1}2)^{\frac {6N-7}{6(6N+5)} + \frac{64}{63}} \geq 2(\tfrac {r_1}2)^{\frac {1}{6} + \frac{64}{63}}
$$
so that, in particular,
$$u, \n d \in L^\infty(Q_{\left[\tfrac {r_1}2\right]^{\frac{149}{126}}}(z_0))\, .$$
\end{remark}

\begin{lemma}\label{corbkoch}
Fix any $\delta \in (\tfrac 12,\tfrac {10}{13})$ and $\bar C\in (0,\infty)$.  There
exist numbers $\e^*=\e^*(\delta, \bar C) \in (0,1)$ and $r^*=r^*(\delta) \in (0,1)$ so small that the
following holds for any fixed open set $\Omega \subseteq \R^3$ and $T\in (0,\infty)$:
\\\\
Set $\Omega_T:=\Omega \times (0,T)$ and suppose
$u,d:\Omega_T  \to \R^3$ and $p: \Omega_T \to \R$ satisfy the assumptions \eqref{enspaces}, \eqref{pspace}, \eqref{divfree},  \eqref{preseq} and \eqref{locenta}.
\\\\
For any $z_0\in \Omega_T$ and $r_1 \in (0,r^*]$ such that $Q_{r_1}(z_0) \subseteq \Omega_T$, if
$$
 \frac 1{r_1^{5 - \delta }} \int_{Q_{r_1}(z_0)} \left(|\n u|^2 + |\n^2 d|^2 + |u|^{\frac {10}3} +|\n d|^{\frac {10}3} + |p|^{\frac {5}3} + |d|^{10}\right) \, dz \leq   \e^*\, ,
$$
then
$$u, \n d \in L^\infty(Q_{\frac {r_0}2}(z_0))
$$
for some small $r_0\in (0,r_1)$.
\end{lemma}
\begin{remark}\label{rmkb}
It will be clear from the proof of Lemma \ref{corbkoch} that its conclusion is true, for example, with
$$r^*:=(\tfrac 14)^{\frac{400}{10-13\delta}} \quad \textrm{and} \quad r_0:=2\left(\tfrac {r_1}2\right)^{\frac{3(10-13\delta)}{400}+1}\geq 2\left(\tfrac {r_1}2\right)^{\frac 74}$$
so that, in particular,
$$u, \n d \in L^\infty (Q_{\left[\frac {r_1}2\right]^{\frac 74}}(z_0))\, .$$
\end{remark}
\noindent
Let us now briefly outline how each theorem (for completeness, we include Theorem \ref{corbliuthm}) follows easily from the corresponding lemma along with  Proposition \ref{covarglem}:
\\\\
{\bf Proofs of Theorems \ref{corbliuthm}, \ref{corbkochbthmb} and \ref{corbkochthm}.} \quad Let us set\footnote{For the proof of Theorem \ref{corbliuthm}, one can in fact remove the term $|d|^{10}$ from $H$, but doing so does not change the conclusion reached by this method of proof.  Moreover, the observation that (in view of \eqref{enspaces}) $H$ satisfies \eqref{partialreggenfb} even if $|d|^{10}$ is included  is crucial to the proofs of Theorems \ref{corbkochbthmb} and \ref{corbkochthm}.}
$$H:=|\n u|^2 + |\n^2 d|^2 + |u|^{\frac {10}3} +|\n d|^{\frac {10}3} + |p|^{\frac {5}3} + |d|^{10}\, .$$
Under the assumptions of Theorem \ref{corbliuthm}, \ref{corbkochbthmb} or \ref{corbkochthm}, it follows (see, e.g., \cite{koch2021}) that
$H$ satisfies \eqref{partialreggenfb} (at least for any bounded $\tilde \Omega$ such that $\mathcal{K} \subset \tilde \Omega_T \subseteq \Omega_T$).
\\\\
Using  Lemma \ref{corbliu} for Theorem \ref{corbliuthm}, Lemma \ref{corbkochb} for Theorem \ref{corbkochbthmb} and Lemma \ref{corbkoch} for Theorem \ref{corbkochthm}, in each setting there exists some $\lambda>0$ ($\lambda=\frac 53 - \gamma$ for some $\gamma \in (0,\tfrac {10}{63})$ in Lemmas \ref{corbliu} and \ref{corbkochb}, $\lambda= 5-\delta$ for some $\delta \in (\tfrac 12,\tfrac {10}{13})$ in Lemma \ref{corbkoch}) as well as some $\e^*>0$ and  $r^*>0$ such that if\footnote{Recall that $Q_{r}(z) \subset Q^*_{r}(z)$ in general.}\\
$$
 \left( \ \frac 1{r_1^{\lambda }} \int_{Q_{r_1}(z_0)} H(z) \, dz
\ \leq\ \right)\ \  \frac 1{r_1^{\lambda }} \int_{Q^*_{r_1}(z_0)} H(z) \, dz
 \leq   \e^*
$$
\ \\\\\\
for some $z_0\in \Omega_T$ and ${r_1} \in (0,r^*]$ such that $Q^*_{r_1}(z_0) \subseteq \Omega_T$,
then $u, \n d \in L^\infty(Q_{\frac {r_0}2}(z_0))$
for some small $r_0\in (0,r_1)$, and hence $z_0 \notin \Sigma_{-}$.  This shows that for any $z_0\in \Sigma_{-}$, we must in particular have
$$\frac 1{r_1^{\lambda}} \int_{Q^*_{r_1}(z_0)} H(z) \, dz
 \geq  \e^*$$
for any such small $r_1>0$, and we may conclude by Proposition \ref{covarglem} with $\mathcal{S}:=\Sigma_{-}$, $\bar r:= r^*$ and\linebreak $c_0:=\e^*$. \hfill $\Box$
\\\\\\
To put these lemmas  (and hence the theorems) into perspective, one should recall the following\footnote{We will see in  Section \ref{epsregsec} below how Lemma \ref{cora} in this form follows easily from the results given in  \cite{koch2021}.} more classical type of  `epsilon-regularity' criteria, due primarily (when \eqref{maineq} also holds) to \cite{caf} (or even \cite{scheffer77}) when $d\equiv 0$ and then extended to general $d\in L^\infty$ in \cite{linliu} and recently by the author in \cite{koch2021} (see Section \ref{epsregsec} below) to more general $d$ when \eqref{maineq} itself need not hold:
\\
\begin{lemma}\label{cora}
Fix any  $\bar C,D\in (0,\infty)$.  There exists $\bar \e=\bar \e (\bar C,D) \in (0,1)$ and, for each $q \in (5,6]$,  there exists $\bar \e_{q} = \bar \e_{q}(\bar C) \in (0,1)$ so small
that the following holds for any fixed open set $\Omega \subseteq \R^3$ and $T\in (0,\infty)$:
\\\\
Set $\Omega_T:=\Omega \times (0,T)$ and suppose
$u,d:\Omega_T  \to \R^3$ and $p: \Omega_T \to \R$ satisfy the assumptions \eqref{enspaces}, \eqref{pspace}, \eqref{divfree},  \eqref{preseq} and \eqref{locenta}.
\\\\
For any $z_0\in \Omega_T$ and $r_0 \in (0,1]$ such that $Q_{r_0}(z_0) \subseteq \Omega_T$, if either
\begin{equation}\label{coraeps}
\!\!\!\!\!\!\!\!\!\!\!\!\!\!\!\frac 1{r_0^{2}} \int_{Q_{r_0}(z_0)} \left(|u|^3 +|\n d|^3 + |p|^{\frac 32}\right) \, dz +  \frac 1{r_0^{2+\frac q2}} \int_{Q_{r_0}(z_0)}|d|^q|\n d|^{3(1-\frac q6)}\, dz \leq  \bar \e_{q} \quad \textrm{for some} \ q \in (5,6]
\end{equation}
or
\begin{equation}\label{coraepsbdd}
d\in L^\infty(\Omega_T) \ \textrm{with} \ \|d\|_{L^\infty(\Omega_T)}\leq D \quad \textrm{and} \quad \frac 1{r_0^{2}} \int_{Q_{r_0}(z_0)} \left(|u|^3 +|\n d|^3 + |p|^{\frac 32}\right) \, dz \leq  \bar \e \, ,
\end{equation}
then
\begin{equation}\label{localboundedness}
u, \n d \in L^\infty(Q_{\frac {r_0}2}(z_0))\, .
\end{equation}
\end{lemma}
\noindent
\ \\\\
Lemmas \ref{corbliu},  \ref{corbkochb}  and  \ref{corbkoch} all rely, in a fundamental way, on some version of Lemma \ref{cora}. For Lemma \ref{corbliu}, Q. Liu (\cite{qliu2018}) used the part for $d\in L^\infty$ in the form given in \cite{linliu} for solutions to \eqref{maineq}, whereas here we rely on \eqref{coraeps} of the more general Lemma \ref{cora}  to prove Lemmas  \ref{corbkochb}  and  \ref{corbkoch}.
\\\\\\
As Lemma \ref{cora} is of a similar form to Lemmas \ref{corbliu},  \ref{corbkochb}  and  \ref{corbkoch}, one may wonder how close one may come to Theorems \ref{corbliuthm}, \ref{corbkochbthmb} and \ref{corbkochthm} using only  Lemma \ref{cora} and Proposition \ref{covarglem} and the method of proof outlined above.  When $d$ is bounded as in Theorem \ref{corbliuthm}, for example, if one uses the epsilon regularity criterion \eqref{coraepsbdd} directly and applies Proposition \ref{covarglem} as in the proof above, one would only reach the conclusion that $\textrm{dim}_{\textrm{pf}}(\Sigma_{-} \cap \mathcal{K}) \leq  2$, as $r_0$ appears in \eqref{coraepsbdd} raised to the power $2$. (One should compare this to the result $\textrm{dim}_{\mathcal{P}}(\Sigma)\leq 2$ for $d\equiv 0$ proved by Scheffer in \cite{scheffer77}.) One can, however, in fact reach $\textrm{dim}_{\textrm{pf}}(\Sigma_{-} \cap \mathcal{K}) \leq  \tfrac 53$ (which should be compared to the improved\footnote{Scheffer's result in \cite{scheffer80} was subsequently improved in \cite{caf} from $\textrm{dim}_{\mathcal{P}}(\Sigma)\leq \tfrac 53$ to $\textrm{dim}_{\mathcal{P}}(\Sigma)\leq 1$, but this required a highly non-trivial application of (some version of) Lemma \ref{cora}.} result $\textrm{dim}_{\mathcal{P}}(\Sigma)\leq \tfrac 53$ for $d\equiv 0$ obtained by  Scheffer in \cite{scheffer80}), corresponding to the lower bound $\gamma=0$ appearing in Lemma \ref{corbliu} and Theorem \ref{corbliuthm}, more or less immediately from \eqref{coraepsbdd} along with H\"older's inequality.  Indeed, one sees easily that Lemma \ref{cora} has the following immediate corollary, which more closely resembles \linebreak Lemmas \ref{corbliu},  \ref{corbkochb}  and  \ref{corbkoch}:
\\\\\\

\pagebreak
\begin{cor}\label{corb}
Fix any  $\bar C,D\in (0,\infty)$.  There exists $\bar \e^*=\bar \e^* (\bar C,D) \in (0,1)$ and, for each $q \in (5,6]$,  there exists $\bar \e_{q}^* = \bar \e_{q}^*(\bar C) \in (0,1)$ so small
that the following holds for any fixed open set $\Omega \subseteq \R^3$ and $T\in (0,\infty)$:
\\\\
Set $\Omega_T:=\Omega \times (0,T)$ and suppose
$u,d:\Omega_T  \to \R^3$ and $p: \Omega_T \to \R$ satisfy the assumptions \eqref{enspaces}, \eqref{pspace}, \eqref{divfree},  \eqref{preseq} and \eqref{locenta}.
\\\\
For any $z_0\in \Omega_T$ and ${r_0} \in (0,1]$ such that $Q_{r_0}(z_0) \subseteq \Omega_T$, if either
\begin{equation}\label{corbeps}
\frac 1{r_0^{\frac 53}} \int_{Q_{r_0}(z_0)} \left(|u|^{\frac {10}3} +|\n d|^{\frac {10}3} + |p|^{\frac {5}3}\right) \, dz\  +\  \frac 1{r_0^{5}} \int_{Q_{r_0}(z_0)}|d|^{10}\, dz \leq \bar   \e^*_q
\quad \textrm{for some}\ q\in (5,6]
\end{equation}
or
\begin{equation}\label{corbepsbdd}
d\in L^\infty(\Omega_T) \ \textrm{with} \ \|d\|_{L^\infty(\Omega_T)}\leq D \quad \textrm{and} \quad \frac 1{r_0^{\frac 53}} \int_{Q_{r_0}(z_0)} \left(|u|^{\frac {10}3} +|\n d|^{\frac {10}3} + |p|^{\frac {5}3}\right) \, dz \leq  \bar  \e^*\, ,
\end{equation}
then
$$u, \n d \in L^\infty(Q_{\frac {r_0}2}(z_0))\, .
$$
\end{cor}
\noindent
{\bf Proof of Corollary \ref{corb}.} \quad Corollary \ref{corb} is a simple consequence of Lemma \ref{cora} along with H\"older's inequality, in view of the convexity estimate
$|d|^q|\n d|^{3(1-\frac q6)} \leq \tfrac q6 |d|^6 + (1-\tfrac q6)|\n d|^3$. \hfill $\Box$
\ \\\\
For $d\in L^\infty(\Omega_T)$, one can now see clearly from the appearance of the power $\tfrac 53$ on $r_0$ in \eqref{corbepsbdd} how one would use Corollary \ref{corb} along with Proposition \ref{covarglem} to obtain the bound $\textrm{dim}_{\textrm{pf}}(\Sigma_{-} \cap \mathcal{K}) \leq  \tfrac 53$.    Therefore the novelty of Q. Liu's Lemma \ref{corbliu}  was that it reached $\tfrac 53-\gamma$ for some $\g >0$ compared to the $\tfrac 53$ (i.e., $\g =0$) which appears in \eqref{corbepsbdd}; more precisely,  Lemma \ref{corbliu} extended this improvement all the way to $\g \approx \tfrac {10}{63}$, so in some sense $\tfrac {10}{63}$ quantifies the improvement\footnote{ However, since the smaller cylinders of radius $\frac {r_0}2$ in
 Lemma \ref{cora} and Corollary \ref{corb} have radii depending linearly on the larger cylinders of radius $r_0$, those results would in fact yield (see e.g. \cite{caf,koch2021} for details) estimates on the dimension of the full singular set $\Sigma$, e.g. $\textrm{dim}_{\textrm{pf}}(\Sigma \cap \mathcal{K}) \leq  \tfrac 53$, where $\Sigma$ is defined similar to $\Sigma_{-}$, but with $Q_r(z_0)$ replaced by the larger, and centered, $Q^*_r(z_0)$.  On the other hand, as pointed out in Remark \ref{remthree}, the smaller radii $\frac{r_0}2$ in Lemma \ref{corbliu} and Lemma \ref{corbkochb} depend non-linearly on $r_1$ ($\frac{r_0}2 \sim (\frac{r_1}2)^{1+\eta}$ for some $\eta>0$) and it therefore seems difficult to conclude from them an estimate on the potentially larger set $\Sigma$ (which would also include points where, roughly speaking, blow-up occurs backwards in time).  Remark \ref{rmkb} shows the situation to be similar in the context of Lemma  \ref{corbkoch}.} of Lemma \ref{corbliu} over Corollary  \ref{corb} (and hence  over Lemma \ref{cora}) in this context. Lemma \ref{corbkochb} similarly obtains this improvement of $\tfrac {10}{63}$, even though one may relax slightly the assumption that $d$ is bounded by assuming \eqref{gsigmadeffina} instead.
\ \\\\
In the case (as in Theorem \ref{corbkochthm}) where (potentially) $d\notin L^\infty(\Omega_T)$, if one instead uses the epsilon regularity criterion \eqref{corbeps} of Corollary \ref{corb} and applies Proposition \ref{covarglem}, then due to the appearance of the larger power $5$ on $r_0$ in \eqref{corbeps} one would only obtain the bound $\textrm{dim}_{\textrm{pf}}(\Sigma_{-} \cap \mathcal{K}) \leq  5$ which would correspond to   $\delta =0$ in the language of Lemma \ref{corbkoch} and Theorem \ref{corbkochthm}.  In that setting, however, the epsilon regularity criterion \eqref{coraeps} (with $q\approx 5$) of Lemma \ref{cora} itself yields (due to the larger power $2+\frac q2$ on $r_0$ in \eqref{coraeps}) the slightly better  bound  of $\textrm{dim}_{\textrm{pf}}(\Sigma_{-} \cap \mathcal{K}) \leq 2+\tfrac 52 = \tfrac 92 = 5-\tfrac 12$ corresponding to  the lower bound $\delta = \tfrac 12$ in Lemma \ref{corbkoch} and Theorem \ref{corbkochthm}. (Note that the criterion \eqref{coraeps} was established in \cite{koch2021} and used there to obtain  the similar estimate  $\textrm{dim}_{\mathcal{P}}(\Sigma)\leq \tfrac 92$, see \cite[Theorem 1]{koch2021}.) On the other hand, applying Proposition \ref{covarglem} to Lemma \ref{corbkoch}  allows one to reach the bound $\textrm{dim}_{\textrm{pf}}(\Sigma_{-} \cap \mathcal{K}) \leq 5-\delta$ even for $\delta \approx \tfrac {10}{13}$ ($\, >\tfrac 12$), and hence Lemma \ref{corbkoch} (which in fact remains true for any $\delta \in (0,\tfrac {10}{13})$) still provides an improvement (in this context), roughly speaking, of $\tfrac {10}{13} - \tfrac 12 = \tfrac {7}{26}$ over Lemma \ref{cora}.

\section{Proof of Proposition \ref{covarglem}}\label{proponepfsec}
In order to prove Proposition \ref{covarglem} (whose proof will certainly already be known to experts -- we include a proof here for clarity and the benefit of the reader),  we first note that the parabolic Vitali covering lemma recalled, for example, in \cite[Lemma 6.1]{caf}, is based on the following fact, which is geometrically clear:
\begin{prop}\label{parainter}
Fix any $\eta \in (0,1)$ and $\beta>0$, and for any $r>0$ and $z:=(x,t)\in \R^3\times \R$ define the parabolic cylinder
$$Q^{\beta,\eta}_{r}(z):=B_r(x)\times I^{\beta,\eta}_{r}(t) \quad \textrm{with} \quad I^{\beta,\eta}_{r}(t):=(t-(1-\eta)\beta r^2,t+\eta \beta r^2)\, .$$
Suppose $A,B >0$ with
$$A\geq \max \left\{2B+1,\sqrt{\frac{B^2}{\min\{\eta,1-\eta\}}+1}\right\}\, .$$
Then for any $z_0,z_1 \in \R^3\times \R$ and $r_0,r_1>0$, if
$$r_1 \leq Br_0 \quad \textrm{and} \quad
Q^{\beta,\eta}_{r_0}(z_0)\cap Q^{\beta,\eta}_{r_1}(z_1) \neq \emptyset
$$
then
$$
Q^{\beta,\eta}_{r_1}(z_1)\subseteq Q^{\beta,\eta}_{Ar_0}(z_0)
\, .$$
\end{prop}
\noindent
The proof is simple, as it is geometrically clear that if even $\overline{B_{r_0}(x_0)}\cap \overline{B_{r_1}(x_1)} \neq \emptyset$ (with the extreme being that the intersection is a single point) and $r_0+2r_1 \leq Ar_0$ then $B_{r_1}(x_1)\subseteq B_{Ar_0}(x_0)$,  and
if even $\overline{I^{\beta,\eta}_{r_0}(t_0)}\cap \overline{I^{\beta,\eta}_{r_1}(t_1)} \neq \emptyset$ (with the extreme being that the intersection is a single point) and both $\eta \beta r_0^2 + \beta r_1^2 \leq \eta \beta  (Ar_0)^2$ as well as
$(1-\eta) \beta r_0^2 + \beta r_1^2 \leq (1-\eta) \beta  (Ar_0)^2$, then $I^{\beta,\eta}_{r_1}(t_1)\subseteq I^{\beta,\eta}_{Ar_0}(t_0)$.  (We omit the remaining details.)
\\\\
In the proof of \cite[Lemma 6.1]{caf}, as they define the cylinders $Q^*:=Q^{1,\frac 18}$, Proposition \ref{parainter} is used with $B:=\frac 32$ (in fact, any $B>1$ would work, with $A$ adjusted appropriately) and with $A:=5$, which is suitable as
$$\max\{2\cdot \tfrac 32 +1,\sqrt{8\cdot (\tfrac 32)^2+1}\} = \sqrt{19} \leq 5\, .$$
However, the proof would similarly go through if we took $Q^*:=Q^{2,\frac 12}$, which are  centered parabolic cylinders of the form (\ref{qstarsdef}),
in which case one may even take $A=4$ (so taking $A=5$ would certainly work as well) as
$$\max\{2\cdot \tfrac 32 +1,\sqrt{2\cdot (\tfrac 32)^2+1}\} = 4\, .$$
We may therefore use \cite[Lemma 6.1]{caf} as stated, but with the parabolic cylinders $Q^{*}_{r}(z)$ defined as in \eqref{qstarsdef} rather than as $Q^{1,\frac 18}_{r}(z)$; let us call such a result $\cite[\textrm{Lemma 6.1}]{caf}_{(\ref{qstarsdef})}$.
\\\\
{\bf Proof of Proposition \ref{covarglem}.} \quad As $\mathcal{K}$ is compact and $\Omega_T$ is open, there exists some $\tilde r \in (0,\bar r]$ sufficiently small that $Q^*_{\tilde r}(z_0)\subseteq \Omega_T$ for any $z_0 \in \mathcal{K}$. Indeed, if not then for any $r\in (0,\bar r]$, there exists $z(r) \in \mathcal{K}$ and $w(r) \notin \Omega_T$ such that $w(r) \in Q^*_{r}(z(r))$.  As $\mathcal{K}$ is (sequentially) compact, there exists a sequence $0< r_k \to 0$ as $k \to \infty$ and $\bar z \in \mathcal{K}$ such that $z(r_k) \to \bar z$ as $k\to \infty$.  By the triangle inequality, this implies as well that $w(r_k) \to \bar z$ as $k\to\infty$, which (as $\bar z \in \mathcal{K} \subset \Omega_T$) contradicts the fact that $\Omega_T$ is open, thus proving the assertion.
\\\\
Fix $\tilde r$ as above. For any $r\in (0,5\tilde r]$, let us consider the cover of the closure ($\overline{\mathcal{S}\cap \mathcal{K}}$) of $\mathcal{S}\cap \mathcal{K}$ by the family of parabolic cylinders centered in $\mathcal{S}\cap \mathcal{K}$ with  radius $\tfrac r5$:
$$\overline{\mathcal{S}\cap \mathcal{K}} \subseteq \bigcup_{z\in \mathcal{S}\cap \mathcal{K}} Q^*_{\frac r5}(z)$$
(as $r>0$ is fixed for all $z$, this  indeed provides a cover of the {\em closure}).
Since $\overline{\mathcal{S}\cap \mathcal{K}}$ is also compact, there exists a finite sub-cover, i.e. there is some $J\in \N$ and a set  of distinct points  $\{ z_j\}_{j=1}^J \subseteq \mathcal{S}\cap \mathcal{K}$ such that
$$\overline{\mathcal{S}\cap \mathcal{K}} \subseteq \bigcup_{j=1}^J Q^*_{\frac r5}( z_j)  \ (\, \subseteq \, \Omega_T \ \textrm{as} \ \tfrac r5 < \tilde r\, )\, .$$
By $\cite[\textrm{Lemma 6.1}]{caf}_{(\ref{qstarsdef})}$, there exists a further sub-set of this finite cover whose members (if there is more than one) are pair-wise disjoint, i.e. there exists $K\in \N$ with $K\leq M$ and a set of strictly increasing sub-indices  $\{j_k\}_{k=1}^K \subseteq \{j\}_{j=1}^M$ such that  $Q^*_{\frac r5}(z_{j_k})\cap Q^*_{\frac r5}(z_{j_\ell}) = \emptyset$ if  $1\leq k < \ell \leq K$, with the additional property that
$$(\, \mathcal{S}\cap \mathcal{K} \, \subseteq \, )\quad \bigcup_{j=1}^M Q^*_{\frac r5}( z_j) \subseteq \bigcup_{k=1}^K Q^*_{r}(z_{j_k})\, ;$$
in particular, the union on the right of $K$ cylinders with radii ($r$) five times those ($\frac r5$) of the original cover constitutes an additional finite open cover.  According to Definition \ref{parfracdimdef}, this implies that
\begin{equation}\label{smalllamnlessa}
N(\mathcal{S}\cap \mathcal{K},r)
\leq K
\, .
\end{equation}
On the other hand, as ${\frac r5} \leq \tilde r \leq \bar r$ and $z_{j_k} \in \mathcal{S}$ for $1\leq k\leq K$, we deduce from \eqref{smalllam} (along with the definition of $\tilde r$) that
$$\frac 1{\left({\frac r5}\right)^{\lambda}}\int_{Q^*_{{\frac r5}}(z_{j_k})} H(z)\, dz \geq   c_0 \quad \textrm{for} \ 1\leq k\leq K\, .
$$
As the cylinders $Q^*_{\frac r5}(z_{j_k})\subseteq \Omega_T$ are disjoint, we then have
$$
K c_0\left({\frac r5}\right)^{\lambda}
=\sum_{k=1}^K c_0\left({\frac r5}\right)^{\lambda}
\leq
\sum_{k=1}^K\int_{ Q^*_{\frac r5}(z_{j_k})} H(z)\, dz
=
\int_{\bigcup_{k=1}^K Q^*_{\frac r5}(z_{j_k})} H(z)\, dz
\leq
\int_{\Omega_T} H(z)\, dz
$$
so that
\begin{equation}\label{smalllamnlessb}
 K
\leq
 c_{\lambda, H}\left({\frac 1r}\right)^{\lambda}
\end{equation}
with
$$
 c_{\lambda, H}:=\frac {5^{\lambda}}{c_0}\int_{\Omega_T}  H(z)\, dz <\infty\, .
$$
According to \eqref{smalllamnlessa} and \eqref{smalllamnlessb} and in view of  Definition \ref{parfracdimdef}, we deduce that
$$\textrm{dim}_{\textrm{pf}}(\mathcal{S}\cap \mathcal{K})=
\limsup_{r\searrow 0} \frac{\log N(\mathcal{S}\cap \mathcal{K},r)}{\log (\tfrac 1r)}
\leq
\limsup_{r\searrow 0} \frac{\log c_{\lambda,H} + \lambda \log (\tfrac 1r)}{\log (\tfrac 1r)} = \lambda
$$
which completes the proof. \hfill $\Box$

\section{Known epsilon-regularity criteria}\label{epsregsec}
\noindent
In this section we show how  Lemma \ref{cora} follows easily from certain results established in \cite{koch2021}.
In the following, for a given $z_0 = (x_0, t_0) \in \rt \times \R$ and $r>0$,  we will adopt the following the notation for the components of a standard parabolic cylinder $Q_r(z_0)$:
\begin{equation}\label{defnqr}
 Q_r(z_0):= B_r(x_0) \times I_r(t_0) \, , \quad \textrm{where} \quad
I_r(t_0):= (t_0 - r^2,t_0) \, .
\end{equation}
Later on, we will also use the shorthand notations (for any set $\mathcal{O}$)
\begin{equation}\label{defnqrlebint}
\left\| U\right\|_{q,\mathcal{O}}:=\left\| U\right\|_{L^q(\mathcal{O})}
\quad \textrm{and} \quad
\left\| V \right\|_{q,s;Q_r(z_0)}:=\left\| V \right\|_{L^s(I_r(t_0);L^q(B_r(x_0)))}
\end{equation}
\ \\\\
The following lemma was proved\footnote{Strictly speaking, in \cite{koch2021}, the first assumption on $d$ in \eqref{f} appears as $d \in L^\infty(I_1(\bar t);L^2(B_1(\bar x)))$ rather than $d \in L^\infty(I_1(\bar t);L^4(B_1(\bar x)))$;  as described in \cite{koch2021}, the $L^4$ assumption is however the most natural one (in view of \eqref{globenineq}, as $F$ is essentially quartic) and of course implies the $L^2$ assumption as the underlying set is bounded.}  in \cite{koch2021}:
\\\\\\

\pagebreak
\begin{lemma}[Lemma 1, \cite{koch2021};  cf. Theorem 2.6 of \cite{linliu} and Proposition 1 of \cite{caf}]\label{thma}
Fix any  $\tilde C,D\in (0,\infty)$. There exists $\tilde \e=\tilde \e (\tilde C,D) \in (0,1)$ and, for each $q \in (5,6]$,  there exists ${\tilde \e_{q} = \tilde \e_{q}(\tilde C) \in (0,1)}$ so small that the following holds for any fixed ${\bar z = (\bar x,\bar t) \in \rt \times \R}$ and $\bar \rho \in (0,1]$:
\\\\
Suppose
$u,d:Q_1(\bar z)  \to \R^3$ and $p: Q_1(\bar z) \to \R$  satisfy (see (\ref{defnqr}))
\begin{equation}\label{f}
\left .\begin{array}{c}
d \in L^\infty(I_1(\bar t);L^4(B_1(\bar x)))\, , \quad
u, \n d \in L^\infty(I_1(\bar t);L^2(B_1(\bar x)))\, , \quad
\nabla u,  \n^2 d  \in L^2(Q_1(\bar z))\\\\
\textrm{and} \quad p\in L^{\frac 32}(Q_1(\bar z))\, ,
\end{array}\right\}
\end{equation}
\begin{equation}\label{udivfree}
\nabla \cdot u =0  \quad \ \textrm{in}\ \  \mathcal{D}'(Q_1(\bar z)) \quad \textrm{and}
\end{equation}
\begin{equation}\label{preseqJ}
 -\Delta p = \nabla \cdot (\nabla^T \cdot [u\otimes u+\n d \odot \n d ]) \quad \ \textrm{in}\ \  \mathcal{D}'(Q_1(\bar z))\, ,
\end{equation}
along with the following local energy inequality:\footnote{See Footnote \ref{fnintdef}.}
\begin{equation}\label{locent}
 \!\!\!\!\!\!\boxed{\begin{array}{l}
 \int_{B_1(\bar x) \times \{t\}} \left(|u|^2 + |\n d|^2\right) \phi\, dx  + \int_{\bar t -1}^t \int_{B_1(\bar x)} \left(|\n u|^2+ |\n^2 d|^2\right) \phi \, dx\, d\tau  \\\\
\qquad \leq  \tilde C\int_{\bar t -1}^t \big\{ \int_{B_1(\bar x)\times \{\tau\}} \left[ \left(|u|^2 + |\n d|^2\right)|\phi_t + \D \phi|  + (|u|^3 + |\n d|^3)|\nabla \phi| + \bar \rho |d|^2|\n d|^2 \phi \right] \, dx \\\\
\qquad \qquad \qquad \qquad \quad +\  \big|\int_{B_1(\bar x) \times \{\tau\}}  pu \cdot \nabla \phi \, dx\big|\ \big\}\, d\tau \\\\
 \textrm{for}\ \textrm{a.e.}\ t\in I_1(\bar t)\  \quad \textrm{and}\quad \forall \ \phi \in \mathcal{C}^\infty_{0}(B_1(\bar x)\times (\bar t -1,\infty)) \ \ \textrm{s.t.}\ \ \phi \geq 0\, .
\end{array}}
\end{equation}
Set\footnote{Note that $E_{3,q}<\infty$ by (\ref{f}) and standard embeddings.}
\begin{equation}\label{etdefn}
E_{3,q}:=\intt{Q_1(\bar z)} (|u|^3+ |\n d|^3+ |p|^{\frac 3 2} +|d|^q|\n d|^{3(1-\frac q6)})\ dz
\end{equation}
and
$$
E_{3}:=\intt{Q_1(\bar z)} (|u|^3+ |\n d|^3+ |p|^{\frac 3 2})\ dz\, .
$$
The following then holds:
\begin{enumerate}
\item If $E_{3,q} \leq \tilde \e_q$ for some $q\in (5,6]$,
then $u, \n d \in L^\infty(Q_{\frac 12}(\bar z))$ with
$$\|u\|_{L^\infty(Q_{1/2}(\bar z))}, \|\n d\|_{L^\infty(Q_{1/2}(\bar z))} \leq {\tilde \e_q}^{2/9}\, .$$

\item (cf. \cite[Theorem 2.6]{linliu})\footnote{An outline of the  proof of this second part was given in the introduction to \cite{koch2021}; see also \cite{linliu} when \eqref{maineq} holds as well.}
If $d\in L^\infty(Q_1(\bar z))$ with $\|d\|_{L^\infty(Q_1(\bar z))} \leq D$ and if $E_{3} \leq \tilde \e$, then
 similarly  $u, \n d \in L^\infty(Q_{\frac 12}(\bar z))$ with
$$\|u\|_{L^\infty(Q_{1/2}(\bar z))}, \|\n d\|_{L^\infty(Q_{1/2}(\bar z))} \leq {\tilde \e}^{2/9}\, .$$
\end{enumerate}

\end{lemma}
\noindent
For completeness, let us now use Lemma \ref{thma} to outline the simple proof (which can essentially be found in \cite{koch2021}) of Lemma \ref{cora}:
\\\\
{\bf Proof of Lemma \ref{cora}.} \quad  If $z_0=(x_0,t_0)$ with $x_0\in \Omega$ and $t_0 \in (0,T)$, setting
\\
\begin{equation}\label{recenteredscalingcor}
\begin{array}{c}
\!\!\!\!\!\!\!\!u_{z_0,{r_0}}(x,t):=r_0u(x_0+{r_0}x, t_0+ r_0^2t)\, , \quad p_{z_0,{r_0}}(x,t):=r_0^2p(x_0+{r_0}x, t_0+ r_0^2t)\\\\
\textrm{and}\quad d_{z_0,{r_0}}(x,t):=d(x_0+{r_0}x, t_0+ r_0^2t)\, ,
\end{array}\end{equation}
\ \\\\\\
it follows (see \cite{koch2021}) from the assumptions of Lemma \ref{cora} that the re-scaled triple $(u_{z_0,{r_0}},d_{z_0,{r_0}},p_{z_0,{r_0}})$  satisfies the conditions of Lemma \ref{thma} with $\bar z:=(0,0)$ and $\bar \rho :=r_0^2$, for some $\tilde C = \tilde C(\bar C)$ depending only on the constant $\bar C$ in \eqref{locenta}.  If we set
$$\bar \e(\bar C,D):=\tilde \e(\tilde C(\bar C),D) \quad \textrm{and} \quad \bar \e_q(\bar C):=\tilde \e_q(\tilde C(\bar C))\, ,$$
then the smallness assumptions in \eqref{coraeps} or \eqref{coraepsbdd} of  Lemma \ref{cora} moreover imply that
$$\frac 1{r_0^{2}} \int_{Q_{r_0}(z_0)} \left(|u|^3 +|\n d|^3 + |p|^{\frac 32}\right) \, dz\  +\  \frac 1{r_0^{2+\frac q2}} \int_{Q_{r_0}(z_0)} \mu \cdot |d|^q|\n d|^{3(1-\frac q6)}\, dz=\qquad \qquad \qquad \qquad $$
$$
\qquad \qquad = \int_{Q_1(0,0)}\left(|u_{z_0,{r_0}}|^3  +|\n d_{z_0,{r_0}}|^3+ |p_{z_0,{r_0}}|^{\frac 32} +\mu \cdot |d_{z_0,{r_0}}|^q|\n d_{z_0,{r_0}}|^{3(1-\frac q6)}\right) \, dz \leq \e_\mu
$$
with
$$\e_\mu
:=(1-\mu)\bar \e(\bar C,D) + \mu \bar \e_q(\bar C)
=(1-\mu)\tilde \e(\tilde C,D) + \mu \tilde \e_q(\tilde C)
$$
for some $\mu \in \{0,1\}$. Lemma \ref{thma} then implies that  $|u_{z_0,{r_0}}|,|\n d_{z_0,{r_0}}|\leq  \e_\mu^{\frac 29}$  on $Q_{\frac 12}(0,0)$, from which
$$|u(y,s)|,|\n d(y,s)| \leq \frac {\e_\mu^{\frac 29}}{r_0} \qquad \textrm{for a.e.}\   (y,s) \in
Q_{\frac {r_0}2}(z_0)$$
(and hence the conclusion of Lemma \ref{cora}) follows in view of \eqref{recenteredscalingcor}.\hfill $\Box$

\section{Proofs of supporting lemmas}\label{lemmasproofssec}
\noindent
In what follows, let us set (recalling \eqref{defnqr})
\begin{equation}\label{athroughgdefa}
\left .\!\!\!\!\!\!\!\!\!\!\!\!\!\!\!\!\!\!\begin{array}{c}
\displaystyle{A_{z_0}(r):= \frac 1r\esssup_{t\in I_r(t_0)}  \int_{B_r(x_0)} \left(|u(t)|^2 + |\n d(t)|^2\right)\, dx \, , \quad }
\\\\
\displaystyle{E_{z_0}(r):= \frac 1r \intt{Q_r(z_0)}\left(|\nabla u|^2+|\n^2 d|^2\right)\, dz \, ,}
\\\\
\displaystyle{C_{z_0}(r):= \frac 1{r^2}\intt{Q_r(z_0)} \left(|u|^3 + |\n d|^3\right)\, dz \, ,\qquad  D_{z_0}(r):= \frac 1{r^2}\intt{Q_r(z_0)}|p|^{3/2}\, dz \, ,}
\\\\
\displaystyle{  \qquad \textrm{and} \qquad  G_{q,z_0}(r):= \frac 1{r^{2+ \frac q2}}  \intt{Q_r(z_0)}|d|^q |\n d|^{3(1-\frac q6)} \, dz\, .}
\end{array} \right\}
\end{equation}
We will make use of the following interpolation-type estimate (see \cite{koch2021}) for the range of the quantities $G_{q,z_0}$ (including $G_{0,z_0}\leq C_{z_0}$), a simple consequence of H\"older's inequality:
\begin{equation}\label{gsiginterpest}
0\leq q \leq \sigma \leq 6 \quad \Longrightarrow \quad G_{q,z_0}(r) \leq G_{\sigma,z_0}^{\frac q\sigma}(r)C_{z_0}^{1-\frac q\sigma}(r) \quad \forall \ r>0\, .
\end{equation}
Let us also set
$$
 \mathcal{E}_{z_0}(r):= \int_{Q_{r}(z_0)} \left( |u|^{\frac {10}3} +|\n d|^{\frac {10}3} + |p|^{\frac {5}3} + |d|^{10}\right) \, dz$$
 and
 $$\mathcal{F}_{z_0}(r):= rE_{z_0}(r)=\int_{Q_{r}(z_0)} \left( |\n u|^{2} +|\n^2 d|^2\right) \, dz \, . $$
Note, in particular, that if we were to assume
\begin{equation}\label{mathcalfsmall}
  \mathcal{F}_{z_0}(R)=  \int_{Q_{R}(z_0)} \left( |\n u|^{2} +|\n^2 d|^2\right) \, dz \leq   \e_* R^{\frac 53 +\mu - \gamma }
\end{equation}
for some $R,\e_* >0$ and $\gamma,\mu \in \R$, then
\begin{equation}\label{mathcalfsmalle}
  E_{z_0}(R)= R^{-1} \int_{Q_{R}(z_0)} \left( |\n u|^{2} +|\n^2 d|^2\right) \, dz \leq   \e_* R^{\frac 23 +\mu - \gamma }\, .
\end{equation}
Lemma \ref{corbkochb} and Lemma \ref{corbkoch} will be consequences of  the following technical lemma:
\begin{lemma}\label{technicallem}
Fix any $\bar C \in (0,\infty)$.  There exits $c=c(\bar C)>1$ and, for any $N\in \N$, there exists $c_N= c_N(\bar C)>1$ such that the following holds for any fixed open set $\Omega \subseteq \R^3$ and $T\in (0,\infty)$:
\\\\
Set  $\Omega_T:=\Omega \times (0,T)$ and suppose
$u,d:\Omega_T  \to \R^3$ and $p: \Omega_T \to \R$ satisfy the assumptions \eqref{enspaces}, \eqref{pspace}, \eqref{divfree},  \eqref{preseq} and \eqref{locenta}.
\\\\
For any $z_0\in \Omega_T$ and $r \in (0,1)$ such that $Q_{2r}(z_0) \subseteq \Omega_T$, if
$$\mathcal{E}_{z_0}(2r),\mathcal{F}_{z_0}(2r) \leq \e_*(2r)^{\frac 53 + \mu -\gamma}$$
for some $\e_* \in (0,1]$ and $\mu,\gamma \geq 0$, then for any $\rho \in (0,2r]$ and $\theta \in (0,\tfrac 14]$, setting
$$\alpha:= \frac{\ln (\tfrac 2\rho)}{\ln (\tfrac 1r)} \geq 1 \quad \textrm{and} \quad \beta:=\frac{\ln (\tfrac 1\theta)}{\ln (\tfrac 1r)}>0$$
(i.e., $\rho=2r^\alpha$ and $\theta = r^\beta$) one has the estimates\footnote{
Note that $2r^{N\beta + \alpha}=\theta^N \rho \leq (\tfrac 14)^N\cdot 2r < 2r$, so that $Q_{2r^{N\beta+\alpha}}(z_0)\subset Q_{2r}(z_0)\subseteq \Omega_T$.
}
\begin{equation}\label{lemfiveestd}
D_{z_0}(2r^{N\beta + \alpha})
\leq c_N\e_*^{\frac 9{10}}(r^{\frac 32 \mathcal{A}}+r^{3 \mathcal{B}}+r^{\frac 32 \mathcal{C}})
\, ,
\end{equation}
\begin{equation}\label{lemfiveestg}
G_{6,z_0}(2r^{N\beta + \alpha}) \leq c\e_*^{\frac 3{5}}r^{3 \tilde {\mathcal{A}}}
\end{equation}
and
\begin{equation}\label{lemfiveestc}
 C_{z_0}(2r^{N\beta + \alpha})\leq c \e_*^{\frac 9{10}}(r^{3 \tilde {\mathcal{B}}}+r^{\frac 32 \tilde {\mathcal{C}}})
 \end{equation}
with
$$\begin{array}{ll}
\mathcal{A}:= \tfrac 23 N\beta  + \tfrac 35 (\mu -\gamma)- [\alpha -1]\, , & \tilde{\mathcal{A}}:= \tfrac 15 (\tfrac 53 + \mu -\gamma) -(N\beta+ \alpha )\, ,\\\\
\mathcal{B}:=(\tfrac N3 -1)\beta  + \tfrac 3{10}\mu - \tfrac 9{20}\gamma + [\alpha-1]\, ,& \tilde{\mathcal{B}}:=   \tfrac 3{10}\mu - \tfrac 9{20}\gamma+(N\beta+ [\alpha -1])\quad \textrm{and}\\\\
\mathcal{C}:=-(N+2)\beta  + \tfrac 13 + \tfrac 45 \mu - \tfrac{19}{20}\gamma- [\alpha -1]\, , \quad & \tilde{\mathcal{C}}:=  \tfrac 13 + \tfrac 45 \mu - \tfrac{19}{20}\gamma-(N\beta+ [\alpha -1])\, ,
\end{array}$$
so that $\tilde{\mathcal{B}}=\mathcal{B} + (\tfrac 23 N + 1)\beta > \mathcal{B} $ and $\tilde{\mathcal{C}}=\mathcal{C}+ 2\beta  > \mathcal{C}$.  In particular (as $r<1$),
\begin{equation}\label{cdsmallcn}
C_{z_0}(2r^{N\beta + \alpha})+D_{z_0}(2r^{N\beta + \alpha})
\stackrel{\eqref{lemfiveestd},\eqref{lemfiveestc}}{\leq} (c+c_N)\e_*^{\frac 9{10}}(r^{\frac 32 \mathcal{A}}+r^{3 \mathcal{B}}+r^{\frac 32 \mathcal{C}})
\end{equation}
and, in view of (\ref{gsiginterpest}),
\begin{equation}\label{gsmallsigma}
G_{q,z_0}(2r^{N\beta + \alpha}) \stackrel{\eqref{lemfiveestc}}{\leq} g_{\sigma}^{\frac q\sigma}[c \e_*^{\frac 9{10}}(r^{3  {\mathcal{B}}}+r^{\frac 32  {\mathcal{C}}})]^{1-\frac q\sigma} \quad
\textrm{as long as}\quad  0\leq q \leq \sigma \leq 6
\end{equation}
(with $g_\sigma \in [0,\infty]$ defined as in (\ref{gsigmadeffina})) and (as $\e_* \leq 1$)
\begin{equation}\label{gsmallsix}
G_{q,z_0}(2r^{N\beta + \alpha})
\stackrel{\eqref{lemfiveestg},\eqref{lemfiveestc}}{\leq} c\e_*^{\frac 3{5}}[r^{3 \tilde {\mathcal{A}}}]^{\frac q6}[r^{3 \tilde {\mathcal{B}}}+r^{\frac 32 \tilde {\mathcal{C}}}]^{1-\frac q6}
\quad
\textrm{as long as}\quad  0\leq q \leq 6\, .
\end{equation}
\end{lemma}
\noindent
We will use \eqref{cdsmallcn}
and \eqref{gsmallsigma}  to prove Lemma \ref{corbkochb}, while we will use \eqref{cdsmallcn} and \eqref{gsmallsix} to
prove Lemma \ref{corbkoch}.
Postponing momentarily the proof of Lemma \ref{technicallem}, let
us use it to prove Lemma \ref{corbkochb} and Lemma \ref{corbkoch}:
\\\\
{\bf Proof of Lemma \ref{corbkochb}.} \quad  Under the assumptions
of Lemma \ref{technicallem}, conclusions \eqref{cdsmallcn}
and \eqref{gsmallsigma} imply that
$$C_{z_0}(r_0)+D_{z_0}(r_0) + G_{q,z_0}(r_0)
\leq (c+c_N)\e_*^{\frac 9{10}}(r^{\frac 32 \mathcal{A}}+r^{3
\mathcal{B}}+r^{\frac 32 \mathcal{C}})+
 g_\sigma^{\frac q\sigma}[c \e_*^{\frac 9{10}}(r^{3  {\mathcal{B}}}+r^{\frac 32  {\mathcal{C}}})]^{1-\frac q\sigma}$$ with
$r_0:=2r^{N\beta + \alpha} < 2r$, as long as $0\leq q \leq \sigma \leq 6$. If we knew that
\begin{equation}\label{positiveconstantsb}
\mathcal{A}, \mathcal{B}, \mathcal{C}\geq 0\, ,
\end{equation}
this along with assumption \eqref{gsigmadeffin} for some $\sigma \in (5,6]$ and $D <\infty$ would imply (as $r,\e_* <1 < c$) that
$$C_{z_0}(r_0)+D_{z_0}(r_0) + G_{q_\sigma,z_0}(r_0)
\leq  \tilde c_{N,\sigma,\bar C,D}\e_*^{\frac 9{10}(1-\frac {q_\sigma} \sigma)} \leq \bar \e_{q_\sigma}
$$
for any fixed $q_\sigma\in (5,\sigma)$, for example $q_\sigma:= \frac{5+\sigma}2$, with $\tilde c_{N,\sigma,\bar C,D}:=(3+2D^{\frac {q_\sigma}\sigma})c(\bar C)+3c_N(\bar C)$, provided that
\begin{equation}\label{choiceofepsstars}
\e_* \leq \left(\frac{\bar \e_{q_\sigma}}{\tilde c_{N,\sigma,\bar C,D}}\right)^{\frac {10}{9(1-\frac {q_\sigma}\sigma)}}=:\e^*=\e^*(N,\sigma, \bar C, D)\, .
\end{equation}
The main conclusion of Lemma \ref{corbkochb} would then follow  from
Lemma \ref{cora} if, for example, $\e_* = \e^*$.  (To prove Lemma \ref{corbkochb}, we will need to take a certain $N=N(\gamma)$, so that $\e^*=\e^*(\gamma,\sigma, \bar C, D)$.)  As we'll soon see, \eqref{positiveconstantsb} can only hold if $\g <\tfrac {10}{63}$, and our goal will  be to show that for any $\bar \g \in (0,1)$, there exists some $N$, $\beta$ and $\alpha$ such that \eqref{positiveconstantsb} holds with $\mu=0$ and $\g =\tfrac {10}{63} \bar \g$.
\\\\
Setting $\mu:=0$ and $\g :=\tfrac {10}{63} \bar \g$ for some $\bar \g >0$, note first that \eqref{positiveconstantsb} would imply that
\begin{equation}\label{bargbetarel}
\bar \g \stackrel{(\mathcal{A},\mathcal{B}\geq 0)}{\leq} 6(N-1)\beta \stackrel{(\mathcal{B},\mathcal{C}\geq 0)}{\leq} \frac{2(N-1)}{2N+9}(3-2\bar \g)
\end{equation}
which (ignoring the intermediate inequality involving $\beta$) implies
\begin{equation}\label{bargbetarelb}
\bar \g \leq 6(N-1)\left[\frac{1}{6N+5}\right]  \leq \frac{2(N-1)}{2N+9}(3-2\bar \g)\, .
\end{equation}
(The two inequalities in \eqref{bargbetarelb} are equivalent, and clearly imply that $\bar \g < 1$ as long as $N\geq 2$.)  If \eqref{bargbetarelb} holds for some $N\in \N \cap [2,\infty)$ and $\bar \g \in (0,1)$, we see easily that  \eqref{bargbetarel} will hold if we  take
$$\beta:= \frac 1{6N+5} \ (\, > \, 0\, ) \qquad \iff \quad N\beta = \frac{1-5\beta}6\, .$$
For such a choice of $\beta$ (and $\mu$ and $\g$), we see that \eqref{positiveconstantsb} says that
\begin{equation}\label{positiveconstantsc}
\tfrac 1{14}\bar \g + \tfrac {23}{18}\beta -\tfrac {1}{18} \stackrel{(\mathcal{B}\geq 0)}{\leq} \alpha -1 \stackrel{(\mathcal{A},\mathcal{C}\geq 0)}{\leq} \tfrac {1}{63} - \tfrac {7}{6}\beta + \min\{\, \tfrac {11}{18}\beta + \tfrac {2}{21}(1-\bar \g),\tfrac {1}{2}\cdot \tfrac {19}{63}(1-\bar \g)\, \}\, .
\end{equation}
Note that
$$\min\{\, \tfrac {11}{18}\beta + \tfrac {2}{21}(1-\bar \g),\tfrac {1}{2}\cdot \tfrac {19}{63}(1-\bar \g)\, \} >0 \quad \textrm{for any} \ \beta >0 \ \textrm{and} \ \bar \g <1\, ,$$
that
$$\tfrac 1{14}\bar \g + \tfrac {23}{18}\beta -\tfrac {1}{18} \leq \tfrac {1}{63} - \tfrac {7}{6}\beta \quad \iff \quad \beta \leq \frac {9(1-\bar \g)}{7\cdot 44}$$
and that
$$\tfrac {1}{63} - \tfrac {7}{6}\beta \geq 0 \quad \iff \quad \beta \leq \frac 2{7\cdot 21}\, ;$$
hence for such $\beta$ (equivalently, $N$), \eqref{positiveconstantsc} will hold with  (for example)
$$\alpha:=1+\left[\tfrac {1}{63} - \tfrac {7}{6}\beta\right] \geq 1\, .$$
Fix now any $\bar \g \in (0,1)$, set
\begin{equation}\label{defnmbarg}
M_{\bar \g}:= \frac 16 \max\left\{\frac{7\cdot 44}{9 (1-\bar \g)},\frac{7\cdot 21}2\right\} \ \ (\ >\,  2 \,  )
\end{equation}
and fix
\begin{equation}\label{defnnrstarbarg}
N=N(\bar \g)\in \N \cap [M_{\bar \g},M_{\bar \g}+1) \quad \textrm{and} \quad r^* = r^*(\bar \g):= (\tfrac 14)^{6M_{\bar \g}+11} \ (\, < \, 1\, )\, .
\end{equation}
(Note that this choice of $N$ implies in particular that $6(1-\bar \g)N > 11 > 6+5\bar \g$ which implies \eqref{bargbetarelb}.)
For such fixed $N=N(\bar \g)$, set also
$$\beta_N:=\frac 1{6N+5}\ (\, > \, 0\, )  \quad \textrm{and} \quad \alpha_N:=\tfrac{64}{63}-\tfrac 76 \beta_N\, ,$$
and note that
$$\frac{1}{6(M_{\bar \g}+1)+5} < \beta_N \leq \frac{1}{6M_{\bar \g}+5} < \min \left\{\frac{9(1-\bar \g)}{7 \cdot 44}, \frac 2{7\cdot 21}\right\}$$
which in turn implies that
$$\alpha_N \geq \tfrac {64}{63} - \tfrac 76\cdot \tfrac{2}{7\cdot 21} =1$$
and hence that \eqref{positiveconstantsc} holds with $\beta:=\beta_N$ and $\alpha:= \alpha_N$.  For any $r>0$, let us now set
$$\theta_{N,r}:= r^{\beta_N} \quad \textrm{and} \quad \rho_{N,r}:= 2r^{\alpha_N}\, .$$
Note that if  $r\in (0,r^*] \subset(0,1)$, then
$$0<\theta_{N,r} \leq (\tfrac 14)^{(6M_{\bar \g}+11)\beta_N} \leq (\tfrac 14)^{\frac{6M_{\bar \g}+11}{6(M_{\bar \g}+1)+5}} = \tfrac 14
\quad \textrm{and} \quad 0<\rho_{N,r} \leq 2r\,  ;
$$
in particular,
$$r_0:=2r^{N\beta_N + \alpha_N}=(\theta_{N,r})^N \rho_{N,r} \leq (\tfrac 14)^N\cdot 2r < 2r\, .$$
Taking $\e^*:=\e^*(N(\tfrac{63}{10}\g),\sigma,\bar C, D)$ and $r^* := r^*(\tfrac{63}{10}\g)$ as in \eqref{choiceofepsstars}, \eqref{defnmbarg} and  \eqref{defnnrstarbarg}  for any $\gamma \in (0,\tfrac {10}{63})$, $\sigma \in (5,6]$ and $\bar C, D \in (0,\infty)$, Lemma \ref{corbkochb} therefore follows from Lemma \ref{technicallem} (with $r:=\frac {r_1}2 \in (0,\frac{r^*}2] \subset (0,1)$, $\e_*:=\e^*$, $\mu:=0$, $\rho:= \rho_{N(\frac{63}{10}\g),\frac {r_1}2}$ and $\theta:=\theta_{N(\frac{63}{10}\g),\frac {r_1}2}$, and taking $q:=\frac{5+\sigma}2$ in \eqref{gsmallsigma}) and  Lemma \ref{cora}. \hfill $\Box$
\ \\\\
{\bf Proof of Lemma \ref{corbkoch}.} \quad Under the assumptions
of Lemma \ref{technicallem}, conclusions \eqref{cdsmallcn} and
\eqref{gsmallsix} imply that
$$C_{z_0}(r_0)+D_{z_0}(r_0) + G_{q,z_0}(r_0)
\leq (c+c_N)\e_*^{\frac 9{10}}(r^{\frac 32 \mathcal{A}}+r^{3
\mathcal{B}}+r^{\frac 32 \mathcal{C}})+
 2c\e_*^{\frac 3{5}}(r^{3 {\mathcal{A}}_q'}+r^{\frac 32
{\mathcal{B}}_q'})$$ with
$$r_0:=2r^{N\beta + \alpha}\, , \quad {\mathcal{A}}_q':=\eta_q \tilde {\mathcal{A}} +
(1-\eta_q)\tilde {\mathcal{B}} \quad \textrm{and} \quad
{\mathcal{B}}_q':=\eta_q \cdot 2\tilde {\mathcal{A}} +
(1-\eta_q)\tilde {\mathcal{C}}\, , \quad \textrm{where} \quad \eta_q:=
\frac q6 \, .$$ If we knew that
\begin{equation}\label{positiveconstantsa}
\mathcal{A}, \mathcal{B}, \mathcal{C}, {\mathcal{A}}_q',
{\mathcal{B}}_q' \geq 0
\end{equation}
for some $q\in (5,6)$, this would imply (as $r,\e_* <1$) that
$$C_{z_0}(r_0)+D_{z_0}(r_0) + G_{q,z_0}(r_0)
\leq  \tilde c_{N,\bar C}\e_*^{\frac 3{5}} \leq \bar \e_q
$$
with $\tilde c_{N,\bar C}:=7c(\bar C)+3c_N(\bar C)$, provided that
\begin{equation}\label{epstardeflemthree}
\e_* \leq \left(\frac{\bar \e_q}{\tilde c_{N,\bar C}}\right)^{\frac 53}=:\e^* = \e^*(N,\bar C,q)\, .
\end{equation}
The main conclusion of Lemma \ref{corbkoch} would then follow from
Lemma \ref{cora} if, for example, $\e_* = \e^*$.  (To prove Lemma \ref{corbkoch}, we will see that one may take a certain $q=q(\delta)$ and $N=3$, so that $\e^* = \e^*(\delta,\bar C)$.)
\ \\\\
We now claim that \eqref{positiveconstantsa} holds for some $q=q(\delta)\in (5,6)$ provided
that $\gamma =0$, $\alpha =1$ and $\mu=\frac {10}3 -\delta$ for
some $\delta <\frac{10}{13}$, i.e. $\mu
> \frac{100}{39}$, for some $\beta>0$ and $N \geq 3$ such that
\begin{equation}\label{betansmallmu}
\beta \in (0,\tfrac \mu 4) \quad \textrm{and} \quad N\beta \in
(0,\tfrac \mu {10})\, .
\end{equation}
It is clear that $\mathcal{A}, \mathcal{B} \geq 0$ if $\gamma =0$,
 $\alpha =1$, $\beta\geq 0$, $\mu \geq 0$ and $N\geq 3$.  Under the same
 assumptions, $\mathcal{C}\geq 0$ provided that
$$(N+2)\beta  \leq \tfrac 13 + \tfrac 45 \mu \, ,$$
and hence $\mathcal{C} \geq 0$ if $\mu >0$ and $\beta$ and $N$ satisfy (for example)
\eqref{betansmallmu}, as then
$$(N+2)\beta  < \tfrac \mu {10} + \tfrac \mu 2 =\tfrac 35 \mu <  \tfrac 13 + \tfrac 45 \mu \, .$$
Now, when $\gamma =0$ and $\alpha =1$ we have
$$\tilde{\mathcal{A}}= -\tfrac 23 + \tfrac 15 \mu  -N\beta\, ,\quad \tilde{\mathcal{B}}=   \tfrac 3{10}\mu +N\beta\quad \textrm{and} \quad \tilde{\mathcal{C}}=  \tfrac 13 + \tfrac 45 \mu -N\beta$$
so that
$$
{\mathcal{A}}_q'=
 [\tfrac 15  \eta_q  +
\tfrac 3{10}(1-\eta_q) ]\mu
-(2\eta_q-1)N\beta  -\tfrac 23\eta_q \geq 0
$$
$$\iff \quad  \mu
 \geq
 \frac{10(2\eta_q-1) }{3-  \eta_q }N\beta
+\frac{20\eta_q}{3(3-  \eta_q) }
= 20\left(\frac{q-3 }{18-  q }\right)N\beta
+\frac{20 q}{54-3q } \approx \frac{40}{13}N\beta + \frac{100}{39}
$$
for $q\approx 5$. Hence for any $\mu > \frac{100}{39}$ and $N\geq 3$, we can choose $q>5$ sufficiently close to $5$ and then $\beta \in (0,\tfrac \mu {10N})$  sufficiently small (depending on $q$ and $N$) to ensure that ${\mathcal{A}}_q'\geq 0$ as well.  For a fixed $\mu > \frac{100}{39}$, one can for example take\footnote{Indeed, setting $\kappa:=\frac{39\mu -100}2$ so that $\mu= \frac {100 + 2\kappa}{39}$, we have
$$\frac{20 q}{54-3q } \leq \frac {100 + \kappa}{39}
\quad
\iff
\quad
q\leq 5 + \frac{13\kappa}{360+\kappa}
$$
which holds for example with $q:= 5 + \frac{\kappa}{360+\kappa}$ as in \eqref{choicenqb}.  Next, to make the sum less than or equal to $\mu= \frac {100 + \kappa}{39}+ \frac { \kappa}{39}$,   we require (taking $N:=3$)
$$
20\left(\frac{q-3 }{18-  q }\right)\cdot 3\beta \leq \frac {\kappa}{39}\quad
\iff
\quad
\beta \leq \frac {\kappa}{60\cdot 39}\left(\frac{18-  q }{q-3 }\right)
$$
which holds for example with $\beta:=\frac{\kappa}{15\cdot 40}$ as in \eqref{choicenqb}; indeed, as $q<6$ and $\kappa >0$, we would then have
$$\beta<\frac{\kappa}{15\cdot 39} \stackrel{(\kappa >0)}{\leq} \frac{\kappa}{60\cdot 39\cdot 3}\left(13 - \frac \kappa{360+\kappa}\right)
=
\frac{\left(18-q\right)\kappa}{60\cdot 39\cdot 3}
\stackrel{(q\leq 6)}{\leq} \frac{\kappa}{60\cdot 39}\left(\frac{18-q}{q-3}\right)\, .
$$
}
\begin{equation}\label{choicenqb}
N:=3\, , \quad q=q_\mu:=5 + \frac{39\mu -100}{720 + (39\mu -100)} \quad \textrm{and} \quad \beta=\beta_\mu:= \frac{39\mu -100}{1200}\, .
\end{equation}
Similarly, we have
$${\mathcal{B}}_q'
=\tfrac 13 -\tfrac 53\eta_q +(\tfrac 45  -\tfrac 25 \eta_q)\mu
-(\eta_q+1)N\beta \geq 0
$$

$$
\iff \quad \mu \geq
\frac 53\cdot \frac{ 5\eta_q  -1}{4  -2 \eta_q} +
 5\cdot \frac{\eta_q+1}{4  -2 \eta_q} N\beta
=\frac 53\cdot \frac{ 5q  -6}{24  -2 q} +
 5\cdot \frac{q+6}{24  -2 q} N\beta
 \approx \frac{95}{42} + \frac{55}{14}N\beta
$$
for $q\approx 5$. As $\frac{95}{42} < \frac{100}{39}$, this should not, in principle, impose additional constraints if $\mu > \frac{100}{39}$ and one can in fact check\footnote{Indeed, again setting $\kappa:=\frac{39\mu -100}2$, we note that
$$\frac 53\cdot \frac{ 5q  -6}{24  -2 q} \leq \frac{100}{39} \ \iff \ q \leq \frac{6\cdot 31}{5\cdot 7}
\quad \textrm{and} \quad  15\cdot \frac{q+6}{24  -2 q} \beta \leq \frac{2\kappa}{39} \ \iff \
 \beta \leq \frac{2\kappa}{15\cdot 39}\cdot \frac{24  -2 q}{q+6}\, .
$$
If we take $N,q,\beta$ again as in (\ref{choicenqb}) and if, for example,  $\frac{100}{39}<\mu < \frac{10}{3}$ (in fact, we could take $\mu$ a bit larger) so that $0< \kappa < 15$, then these will both hold (and hence ${\mathcal{B}}_q'\geq 0$, as $\mu = \frac{100}{39}+\frac{2\kappa}{39}$) as then
$$q= 5 + \frac{\kappa}{360+\kappa} < 5 + \frac{15}{360}  < 5 + \frac{1}{7}  <\frac{6\cdot 31}{5\cdot 7}
$$
and, as in particular $q<6$,
$$\beta=\frac{\kappa}{15\cdot 40}<
\frac{2\kappa}{15\cdot 39} \cdot \frac{24  -2 q}{q+6}\, .$$}
that ${\mathcal{B}}_q'\geq 0$ as well  for $N$, $q$ and $\beta$ as in \eqref{choicenqb}.
\\\\
Setting $\mu(\delta):=\frac{10}3 -\delta$ and taking $\e^*:=\e^*(3,\bar C,q_{\mu(\delta)})=\e^*(3,\bar C,5 + \frac{10-13\delta}{250-13\delta})$  and ${r^*:=(\tfrac 14)^{\frac 1{\beta_{\mu(\delta)}}}=(\tfrac 14)^{\frac{400}{10-13\delta}}}$   as in \eqref{epstardeflemthree} and \eqref{choicenqb} for any $\delta \in (\tfrac 12,\tfrac {10}{13})$ (or even smaller) and $\bar C\in (0,\infty)$,  Lemma \ref{corbkoch} now follows from Lemma \ref{technicallem} (with $N:=3$,  $r:=\frac{r_1}2$, $\e_*:=\e^*$ $\mu:=
\frac{10}3-\delta$, $\g:=0$, $\rho:=r_1 = 2r$  and $\theta:= \left(\tfrac{r_1}2\right)^{\frac{10-13\delta}{400}}$ so that $\alpha=1$ and $\beta =
\frac{10-13\delta}{400}$, and taking $q:=5 + \frac{10-13\delta}{250-13\delta}$ in \eqref{gsmallsix}) and Lemma \ref{cora}. \hfill $\Box$
\\\\
To prove Lemma \ref{technicallem}, we will rely crucially on
the following proposition (which is only slightly different from the corresponding result\footnote{One need only make a slight adjustment in the proof of Proposition \ref{propkocha} (and take $\mu=0$) to essentially recover Liu's result, at the point where the term  $\left\| |d|^2 |\n d|^2 \right\|_{1,Q_\rho(z_0)}$ (coming from the local energy inequality) appears. As $|d|^2$ appears as a multiplicative factor, one can simply ignore it at the expense of allowing the constant $C$ in Proposition \ref{propkocha} to depend on $\|d\|_{L^\infty}$; one would then no longer need to include the term $|d|^{10}$ in \eqref{mathcalesmall}.} in \cite{qliu2018}) which is a consequence of the local
energy inequality.  (The reader should recall the notation (\ref{defnqrlebint}) for Lebesgue norms which we will use in all of what follows.)

\begin{prop}\label{propkocha}
There exists $C >0$
such that the following holds:
\\\\
Fix an open set $\Omega \subseteq \R^3$ and $\bar C,T\in (0,\infty)$, set  $\Omega_T:=\Omega \times (0,T)$ and suppose
$u,d:\Omega_T  \to \R^3$ and $p: \Omega_T \to \R$ satisfy the assumptions \eqref{enspaces}, \eqref{pspace} and \eqref{locenta}.
\\\\
Suppose
\begin{equation}\label{mathcalesmall}
  \mathcal{E}_{z_0}(2r)=  \int_{Q_{2r}(z_0)} \left( |u|^{\frac {10}3} +|\n d|^{\frac {10}3} + |p|^{\frac {5}3} + |d|^{10}\right) \, dz \leq   \e_* (2r)^{\frac 53 +\mu - \gamma }
\end{equation}
for some $z_0\in \Omega_T$ and $r \in (0,1]$ such that $Q_{2r}(z_0) \subseteq \Omega_T$ and some $\mu, \gamma \geq 0$ (one may assume $\mu \cdot \gamma =0$) and $\e^* \in (0,1]$.
Then\footnote{Similarly (and for the same reason), $E_{z_0}(r)=r^{-1}\left\| |\n u|^2 + |\n^2 d|^2 \right\|_{2,Q_r(z_0)} \leq C\bar C  \e_*^{\frac 35} r^{\frac 35 \mu - \frac 9{10}\g}$, but this will not help us.} (see (\ref{athroughgdefa}))
\begin{equation}\label{propaest}
A_{z_0}(r)=r^{-1}\left\| |u|^2 + |\n d|^2 \right\|_{1,\infty;Q_r(z_0)} \leq  C\bar C   \e_*^{\frac 35} r^{\frac 35 \mu - \frac 9{10}\g}\, .
\end{equation}
\end{prop}
\noindent
{\bf Proof of Proposition \ref{propkocha}.} \quad Using the backwards heat kernel and a suitable cut-off function, it is not hard to see (see, for example, \cite{koch2021,qliu2018} for more details) that for any $z_0 \in \Omega_T$ and $0 < r \leq \frac \rho 2 \leq 1$, one can construct a test function $0 \leq \phi \in \mathcal{C}^\infty_0(Q_\rho(z_0))$ with the following properties:

\begin{equation}\label{phipropb}
 \frac 1{r^3} \lesssim \phi \qquad \textrm{on} \quad Q_r(z_0)
\end{equation}
and
\begin{equation}\label{phipropa}
\phi \lesssim \frac 1{r^3}\, , \qquad |\n \phi| \lesssim \left(\frac \rho r + 1\right)\frac 1{r^4} \qquad \textrm{and} \qquad |\phi_t + \Delta \phi|\lesssim \frac 1{\rho^5} \qquad \textrm{on} \quad Q_\rho(z_0)\, .
\end{equation}
Applying \eqref{locenta} for such a $\phi$, we have (with the constant in the second inequality being $\bar C$)
$$A_{z_0}(r)=
r^{-1}\left\| |u|^2 + |\n d|^2 \right\|_{1,\infty;Q_r(z_0)}
\stackrel{\eqref{phipropb}}{\lesssim}
r^{2}\left\| (|u|^2 + |\n d|^2)\phi \right\|_{1,\infty;Q_\rho(z_0)}
$$
$$
\stackrel{\eqref{locenta}}{\lesssim}
r^2\left\|\left(|u|^2 + |\n d|^2\right)|\phi_t + \D \phi|  +\left(|u|^2 + |\n d|^2 + |p|\right)|u||\n \phi|  + (1+|d|^2) |\n d|^2\phi \right\|_{1,Q_\rho(z_0)}
$$
$$
\stackrel{\eqref{phipropa}}{\lesssim}
\left(\frac {r^2}{\rho^{5}}+ \frac 1{r}\right)
\left\||u|^2 + |\n d|^2\right\|_{1,Q_\rho(z_0)}+
\left(\frac {\rho}{r}+ 1\right)\frac 1{r^2}
\left\||u|^3 + |\n d|^3 + |p|^{\frac 32}\right\|_{1,Q_\rho(z_0)}+ \frac 1r
\left\| |d|^2 |\n d|^2 \right\|_{1,Q_\rho(z_0)}
$$
$$
\!\!\!\!\!\!\!\!\!\!\lesssim
\left(\frac {r^2}{\rho^{5}}+ \frac 1{r}\right)\rho^2
\left\||u|^2 + |\n d|^2\right\|_{\frac 53,Q_\rho(z_0)}+
\left(\frac {\rho}{r}+ 1\right)\frac {\rho^{\frac 12}}{r^2}
\left\||u|^3 + |\n d|^3 + |p|^{\frac 32}\right\|_{\frac {10}9,Q_\rho(z_0)}+ \frac \rho r
\left\| |d|^2  \right\|_{5,Q_\rho(z_0)}
\left\| |\n d|^2 \right\|_{\frac 53,Q_\rho(z_0)}
$$
$$
\!\!\!\!\!\!\!\!\!\!\lesssim
\left(\left(\frac {r}{\rho}\right)^3+ \left(\frac {\rho}{r}\right)^2 r^2\right)
\left[r^{-\frac 53}
\mathcal{E}_{z_0}(\rho)\right]^{\frac 35}+
\left(\frac {\rho}{r}+ 1\right)\left(\frac {\rho}{r} \right)^{\frac 12}\left[r^{-\frac 53}
\mathcal{E}_{z_0}(\rho)\right]^{\frac 9{10}} + \left(\frac \rho r \right)
r^{\frac 43}\left[r^{-\frac 53}
\mathcal{E}_{z_0}(\rho)\right]^{\frac 45}\, .$$
Therefore, taking $\rho:=2r$ so that, in particular,
$$ \frac \rho r, \frac r\rho , r, \e_*\leq 2\, ,$$
assumption \eqref{mathcalesmall} implies that
$$A_{z_0}(r)\lesssim
\left[\e_* r^{\mu -\gamma}\right]^{\frac 35}+
\left[\e_* r^{\mu -\gamma}\right]^{\frac 9{10}} +
\left[\e_* r^{\mu -\gamma}\right]^{\frac 45}
=\left(
r^{\frac 3{10}\gamma} +\e_*^{\frac 3{10}}r^{\frac 3{10}\mu}+\e_*^{\frac 15}r^{\frac 15\mu +\frac 1{10}\gamma}\right)
\e_*^{\frac 35}r^{\frac 35\mu -\frac 9{10}\gamma} \lesssim \e_*^{\frac 35}r^{\frac 35\mu -\frac 9{10}\gamma}$$
(as $\mu,\gamma \geq 0$).  If we keep track of the constants involved, the conclusion of the proposition follows.  \hfill $\Box$
\\\\
As in \cite{qliu2018}, we will also need the following two propositions, the first of which is a consequence of the Sobolev embeddings and Poincar\'e inequality while the second is a consequence of the elliptic theory; both can be found in the literature (for example, as indicated):
\begin{prop}[Interpolation inequality; see, e.g., \cite{choeyang2015}, Lemma 2]\label{sobprop}
There exists a constant $C>0$ such that, for any $0<  r' \leq r <\infty$ and any measurable function $U:Q_r(z_0) \to \R^3$, the estimate
$$\|U\|_{3,Q_{r'}(z_0)} \leq C \left(
r^{\frac 16}\|U\|_{2,\infty;Q_r(z_0)}^{\frac 12}\|\n U\|_{2,Q_r(z_0)}^{\frac 12}
+ r^{\frac 53}r^{-\frac 32}\|U\|_{2,\infty;Q_r(z_0)}\right)
$$
holds provided the right-hand side is well-defined.  In particular, for any $r>0$ and $\eta \in (0,1]$, one has an estimate of the form (see (\ref{athroughgdefa}))
\begin{equation}\label{caeest}
C_{z_0}(\eta r) \lesssim \eta^{3}A_{z_0}^{\frac 32}(r) + \eta^{-\frac 32}A_{z_0}^{\frac 34}(r)E_{z_0}^{\frac 34}(r) \quad \textrm{for any}\ r>0,\ \eta \in (0,1], \ z_0\in \R^{3+1}
\end{equation}
provided that the right-hand side  is well-defined.
\end{prop}

\begin{prop}[Interior elliptic estimate, see\footnote{\cite{gkt2007} states the result for any $\theta\in (0,\frac 12]$, but the author believes this to be a typographical error.}, e.g., \cite{gkt2007}, Lemma 3.4]\label{ellipprop}
For any $q\in (1,\infty)$ and $n\in \N$, there exists a constant $C_{q,n}>0$ such that if $U$ is a weak solution to $-\Delta U = \nabla \cdot (\nabla^T \cdot F)$ in $B_R(x_0)\subset \R^n$ for some $R>0$ and $x_0\in \R^n$, then
$$\|U\|_{q,B_{\theta R}(x_0)} \leq C_{q,n} \left(\|F\|_{q,B_{R}(x_0)}
+ \theta^{\frac nq}\|U\|_{q,B_{\frac R2}(x_0)}\right) \quad \textrm{for any} \ \theta \in (0,\tfrac 14] \, ,
$$
provided the right-hand side is well-defined.  In particular (for $n=3$ and $q=\tfrac 32$),
\begin{equation}\label{ddcest}
D_{z_0}(\theta R) \lesssim
\theta D_{z_0}(R) + \theta^{-2}C_{z_0}(R)\quad \textrm{for any} \ \theta \in (0,\tfrac 14]
\end{equation}
(see (\ref{athroughgdefa})) provided $p$ satisfies the pressure equation \eqref{preseq} and $Q_R(z_0)\subseteq \Omega_T$.
\end{prop}
\noindent
\ \\
As in \cite{qliu2018}, let us now use these propositions  to prove Lemma \ref{technicallem}.
\ \\\\
{\bf Proof of Lemma \ref{technicallem}.}\quad  Under the assumptions of Lemma \ref{technicallem}, taking $\eta:=\theta^j \tfrac \rho r$ in \eqref{caeest} for any $j\in \N$, Proposition \ref{propkocha} and Proposition \ref{sobprop} imply that
\begin{equation}\label{cthjrhoest}
\left .\begin{array}{rcl}
C_{z_0}(\theta^j \rho) & \stackrel{\eqref{caeest}}{\lesssim} & (\theta^j \tfrac \rho r)^{3}A_{z_0}^{\frac 32}(r) + (\theta^j \tfrac \rho r)^{-\frac 32}A_{z_0}^{\frac 34}(r)E_{z_0}^{\frac 34}(r)
\\\\
 & \stackrel{\eqref{propaest},\eqref{mathcalfsmalle}}{\lesssim_{\bar C}} & (\theta^j \tfrac \rho r)^{3}\left[\e_*^{\frac 35} r^{\frac 35 \mu - \frac 9{10}\g}\right]^{\frac 32} + (\theta^j \tfrac \rho r)^{-\frac 32}\left[\e_*^{\frac 35} r^{\frac 35 \mu - \frac 9{10}\g}\right]^{\frac 34}\left[\e_* r^{\frac 23 +\mu - \gamma }\right]^{\frac 34}
\\\\
 & \stackrel{(\e_* \leq 1)}{\lesssim} & \e_*^{\frac 9{10}}\left[(\theta^j \rho )^{3} r^{\frac 9{10} \mu - \frac {27}{20}\g -3} + (\theta^j  \rho )^{-\frac 32}r^{2+\frac 65 \mu - \frac {3}{2}\cdot \frac {19}{20}\g}
 \right]
\end{array}
\right\}\end{equation}
which, when $j=N$, implies \eqref{lemfiveestc} for some $c=c(\bar C)>1$.  Moreover, noting that
\begin{equation}\label{drhoest}
\!\!\!\!\!D_{z_0}(\rho)\lesssim \rho^{-\frac 32} \left\||p|^{\frac 53}\right\|_{1,Q_\rho(z_0)}^{\frac 9{10}}
\leq
\rho^{-\frac 32} \left\||p|^{\frac 53}\right\|_{1,Q_{2r}(z_0)}^{\frac 9{10}}
\leq
\rho^{-\frac 32} \left[\mathcal{E}_{z_0}(2r) \right]^{\frac 9{10}}
\leq
\rho^{-\frac 32} \left[\e_*(2r)^{\frac 53 + \mu -\gamma} \right]^{\frac 9{10}}\, ,
\end{equation}
repeated applications of Proposition \ref{ellipprop} (with $R:=\theta^{k-1} \rho$ for $k\leq N$) along with \eqref{cthjrhoest} imply that
$$
\begin{array}{rcl}
\!\!\!\!\!\!\!\!\!\!D_{z_0}(\theta^N \rho) & \stackrel{\eqref{ddcest}}{\lesssim_N} &
\displaystyle{\theta^N D_{z_0}(\rho) + \theta^{N-3}\sum_{j=0}^{N-1}\theta^{-j}C_{z_0}(\theta^j \rho)}
\\\\
 & \stackrel{{\tiny
 \begin{array}{c}\eqref{cthjrhoest},\\ \eqref{drhoest} \end{array}
} }{\lesssim_{\bar C}} &
\displaystyle{\e_*^{\frac 9{10}} \left(\theta^N \rho^{-\frac 32} \left[\e_*(2r)^{\frac 53 + \mu -\gamma} \right]^{\frac 9{10}} +
\theta^{N-3}\sum_{j=0}^{N-1}\left[\theta^{2j}\rho ^{3} r^{\frac 9{10} \mu - \frac {27}{20}\g -3} + \theta^{-\frac 52 j}\rho^{-\frac 32}r^{2+\frac 65 \mu - \frac {3}{2}\cdot \frac {19}{20}\g}
 \right]\right)}
\\\\
 & \stackrel{(\theta < 1)}{\leq} &
\displaystyle{\e_*^{\frac 9{10}} \left(\theta^N \rho^{-\frac 32} \left[\e_*(2r)^{\frac 53 + \mu -\gamma} \right]^{\frac 9{10}} +
N\theta^{N-3}\left[\rho ^{3} r^{\frac 9{10} \mu - \frac {27}{20}\g -3} + \theta^{-\frac 52 N}\rho^{-\frac 32}r^{2+\frac 65 \mu - \frac {3}{2}\cdot \frac {19}{20}\g}
 \right]\right)}
\end{array}
$$
which\footnote{Here, as in \cite{qliu2018}, we conclude  very crudely using $\ell^\infty \hookrightarrow \ell^1_{\textrm{loc}}$; that is, we do not make use of any convergent geometric series.} implies \eqref{lemfiveestd} for some $c_N=c_N(\bar C)>1$.  Finally, setting $r_0:=\theta^N \rho \leq (\tfrac 14)^N (2r) < 2r$, we have
$$
G_{6,z_0}(r_0)
\lesssim
r_0^{-3} \left\||d|^{10}\right\|_{1,Q_{r_0}(z_0)}^{\frac 3{5}}
\leq
r_0^{-3} \left\||d|^{10}\right\|_{1,Q_{2r}(z_0)}^{\frac 3{5}}
\leq
r_0^{-3} \left[\mathcal{E}_{z_0}(2r) \right]^{\frac 3{5}}
\leq
r_0^{-3} \left[\e_*(2r)^{\frac 53 + \mu -\gamma} \right]^{\frac 3{5}}
$$
which implies \eqref{lemfiveestg} for some $c>1$ (which in fact does not depend on $\bar C$). \hfill $\Box$


\begin{thebibliography}{GKT07}

\bibitem[CKN82]{caf}
L.~Caffarelli, R.~Kohn, and L.~Nirenberg.
\newblock Partial regularity of suitable weak solutions of the
  {N}avier-{S}tokes equations.
\newblock {\em Comm. Pure Appl. Math.}, 35(6):771--831, 1982.

\bibitem[CY15]{choeyang2015}
Hi~Jun Choe and Minsuk Yang.
\newblock Hausdorff measure of the singular set in the incompressible
  magnetohydrodynamic equations.
\newblock {\em Comm. Math. Phys.}, 336(1):171--198, 2015.

\bibitem[GKT07]{gkt2007}
Stephen Gustafson, Kyungkeun Kang, and Tai-Peng Tsai.
\newblock Interior regularity criteria for suitable weak solutions of the
  {N}avier-{S}tokes equations.
\newblock {\em Comm. Math. Phys.}, 273(1):161--176, 2007.

\bibitem[Koc21]{koch2021}
Gabriel~S. Koch.
\newblock Partial regularity for {N}avier-{S}tokes and liquid crystals
  inequalities without maximum principle.
\newblock {\em arXiv:2001.04098 {\em (to appear in} Analysis \& PDE{\em)}},
  2021.

\bibitem[KP12]{kukpei2012}
Igor Kukavica and Yuan Pei.
\newblock An estimate on the parabolic fractal dimension of the singular set
  for solutions of the {N}avier-{S}tokes system.
\newblock {\em Nonlinearity}, 25(9):2775--2783, 2012.


\bibitem[Liu18]{qliu2018}
Qiao Liu.
\newblock Dimension of singularities to the 3d simplified nematic liquid
  crystal flows.
\newblock {\em Nonlinear Anal. Real World Appl.}, 44:246--259, 2018.

\bibitem[Liu21]{qliu2021}
Qiao Liu.
\newblock Partial regularity and the {M}inkowski dimension of singular points
  for suitable weak solutions to the 3{D} simplified {E}ricksen-{L}eslie
  system.
\newblock {\em Discrete Contin. Dyn. Syst.}, 41(9):4397--4419, 2021.

\bibitem[LL95]{linliu95}
Fang-Hua Lin and Chun Liu.
\newblock Nonparabolic dissipative systems modeling the flow of liquid
  crystals.
\newblock {\em Comm. Pure Appl. Math.}, 48(5):501--537, 1995.

\bibitem[LL96]{linliu}
Fang-Hua Lin and Chun Liu.
\newblock Partial regularity of the dynamic system modeling the flow of liquid
  crystals.
\newblock {\em Discrete Contin. Dynam. Systems}, 2(1):1--22, 1996.

\bibitem[Sch77]{scheffer77}
Vladimir Scheffer.
\newblock Hausdorff measure and the {N}avier-{S}tokes equations.
\newblock {\em Comm. Math. Phys.}, 55(2):97--112, 1977.

\bibitem[Sch80]{scheffer80}
Vladimir Scheffer.
\newblock The {N}avier-{S}tokes equations on a bounded domain.
\newblock {\em Comm. Math. Phys.}, 73(1):1--42, 1980.

\bibitem[Sch85]{scheffer3}
Vladimir Scheffer.
\newblock A solution to the {N}avier-{S}tokes inequality with an internal
  singularity.
\newblock {\em Comm. Math. Phys.}, 101(1):47--85, 1985.

\end{thebibliography}
\end{document}